\documentclass[oneside,english]{amsart}
\usepackage{lmodern}
\usepackage[T1]{fontenc}
\usepackage[latin9]{inputenc}
\usepackage{color}
\usepackage{enumitem}
\usepackage{amsthm}
\usepackage{amssymb}
\usepackage{wasysym}
\PassOptionsToPackage{normalem}{ulem}
\usepackage{ulem}

\makeatletter
\numberwithin{equation}{section}
\numberwithin{figure}{section}
\theoremstyle{plain}
\newtheorem{thm}{\protect\theoremname}[section]
  \theoremstyle{definition}
  \newtheorem{defn}[thm]{\protect\definitionname}
  \theoremstyle{plain}
  \newtheorem{conjecture}[thm]{\protect\conjecturename}
  \theoremstyle{plain}
  \newtheorem{fact}[thm]{\protect\factname}
  \theoremstyle{remark}
  \newtheorem{claim}[thm]{\protect\claimname}
  \theoremstyle{plain}
  \newtheorem{cor}[thm]{\protect\corollaryname}
 \theoremstyle{plain}
 \newtheorem{obs}[thm]{Observation}
  \theoremstyle{plain}
  \newtheorem{prop}[thm]{\protect\propositionname}
 \newlist{casenv}{enumerate}{4}
 \setlist[casenv]{leftmargin=*,align=left,widest={iiii}}
 \setlist[casenv,1]{label={{\itshape\ \casename} \arabic*.},ref=\arabic*}
 \setlist[casenv,2]{label={{\itshape\ \casename} \roman*.},ref=\roman*}
 \setlist[casenv,3]{label={{\itshape\ \casename\ \alph*.}},ref=\alph*}
 \setlist[casenv,4]{label={{\itshape\ \casename} \arabic*.},ref=\arabic*}
  \theoremstyle{remark}
  \newtheorem{notation}[thm]{\protect\notationname}
  \theoremstyle{remark}
  \newtheorem{rem}[thm]{\protect\remarkname}
  \theoremstyle{remark}
  \newtheorem*{claim*}{\protect\claimname}
 \theoremstyle{plain}
 \newtheorem{const}[thm]{Construction}
  \theoremstyle{definition}
  \newtheorem{problem}[thm]{\protect\problemname}
 \theoremstyle{definition}
 \newtheorem*{defn*}{\protect\definitionname}
  \theoremstyle{plain}
  \newtheorem*{thm*}{\protect\theoremname}
  \theoremstyle{definition}
  \newtheorem{example}[thm]{\protect\examplename}

\usepackage{amsmath}

\usepackage{euler}

\usepackage{url}
\theoremstyle{definition}

\usepackage{a4wide}
\usepackage{amssymb}
\usepackage{url}
\makeatother

\usepackage{babel}
  \providecommand{\claimname}{Claim}
  \providecommand{\conjecturename}{Conjecture}
  \providecommand{\corollaryname}{Corollary}
  \providecommand{\definitionname}{Definition}
  \providecommand{\examplename}{Example}
  \providecommand{\factname}{Fact}
  \providecommand{\notationname}{Notation}
  \providecommand{\problemname}{Problem}
  \providecommand{\propositionname}{Proposition}
  \providecommand{\remarkname}{Remark}
  \providecommand{\theoremname}{Theorem}
 \providecommand{\casename}{Case}
\providecommand{\theoremname}{Theorem}

\begin{document}
\global\long\def\p{\mathbf{p}}
\global\long\def\q{\mathbf{q}}
\global\long\def\C{\mathfrak{C}}
\global\long\def\SS{\mathcal{P}}
 \global\long\def\pr{\operatorname{pr}}
\global\long\def\image{\operatorname{im}}
\global\long\def\otp{\operatorname{otp}}
\global\long\def\dec{\operatorname{dec}}
\global\long\def\suc{\operatorname{suc}}
\global\long\def\pre{\operatorname{pre}}
\global\long\def\qe{\operatorname{qf}}
 \global\long\def\ind{\operatorname{ind}}
\global\long\def\Nind{\operatorname{Nind}}
\global\long\def\lev{\operatorname{lev}}
\global\long\def\Suc{\operatorname{Suc}}
\global\long\def\HNind{\operatorname{HNind}}
\global\long\def\minb{{\lim}}
\global\long\def\concat{\frown}
\global\long\def\cl{\operatorname{cl}}
\global\long\def\tp{\operatorname{tp}}
\global\long\def\id{\operatorname{id}}
\global\long\def\cons{\left(\star\right)}
\global\long\def\qf{\operatorname{qf}}
\global\long\def\ai{\operatorname{ai}}
\global\long\def\dtp{\operatorname{dtp}}
\global\long\def\acl{\operatorname{acl}}
\global\long\def\nb{\operatorname{nb}}
\global\long\def\limb{{\lim}}
\global\long\def\leftexp#1#2{{\vphantom{#2}}^{#1}{#2}}
\global\long\def\intr{\operatorname{interval}}
\global\long\def\atom{\emph{at}}
\global\long\def\I{\mathfrak{I}}
\global\long\def\uf{\operatorname{uf}}
\global\long\def\ded{\operatorname{ded}}
\global\long\def\Ded{\operatorname{Ded}}
\global\long\def\Df{\operatorname{Df}}
\global\long\def\Th{\operatorname{Th}}
\global\long\def\eq{\operatorname{eq}}
\global\long\def\Aut{\operatorname{Aut}}
\global\long\def\ac{ac}
\global\long\def\DfOne{\operatorname{df}_{\operatorname{iso}}}
\global\long\def\modp#1{\pmod#1}
\global\long\def\sequence#1#2{\left\langle #1\left|\,#2\right.\right\rangle }
\global\long\def\set#1#2{\left\{  #1\left|\,#2\right.\right\}  }
\global\long\def\Diag{\operatorname{Diag}}
\global\long\def\Nn{\mathbb{N}}
\global\long\def\mathrela#1{\mathrel{#1}}
\global\long\def\twiddle{\mathord{\sim}}
\global\long\def\mathordi#1{\mathord{#1}}
\global\long\def\Qq{\mathbb{Q}}
\global\long\def\dense{\operatorname{dense}}
 \global\long\def\cof{\operatorname{cof}}
\global\long\def\tr{\operatorname{tr}}
\global\long\def\treeexp#1#2{#1^{\left\langle #2\right\rangle _{\tr}}}

\title{Examples in dependent theories}

\author{Itay Kaplan and Saharon Shelah}

\thanks{The first author's research was partially supported by the SFB 878
grant.}

\thanks{The first author was partially supported by SFB grant 878. The second
author would like to thank the Israel Science Foundation for partial
support of this research (Grants nos. 710/07 and 1053/11). No. 946
on the second author's list of publications.}

\address{Itay Kaplan \\
 Institut f\"ur mathematische Logik und Grundlagenforschung \\
 Fachbereich Mathematik und Informatik \\
 Universit\"at M\"unster \\
 Einsteinstra{\ss}e 62 \\
 48149 M\"unster \\
 Germany }

\email{itay.kaplan@uni-muenster.de}

\urladdr{https://sites.google.com/site/itay80/ }

\address{Saharon Shelah\\
The Hebrew University of Jerusalem\\
Einstein Institute of Mathematics \\
Edmond J. Safra Campus, Givat Ram\\
Jerusalem 91904, Israel}

\address{Saharon Shelah \\
Department of Mathematics\\
Hill Center-Busch Campus\\
Rutgers, The State University of New Jersey\\
110 Frelinghuysen Road\\
Piscataway, NJ 08854-8019 USA}

\email{shelah@math.huji.ac.il}

\urladdr{http://shelah.logic.at/}

\begin{abstract}
In the first part we show a counterexample to a conjecture by Shelah
regarding the existence of indiscernible sequences in dependent theories
(up to the first inaccessible cardinal). In the second part we discuss
generic pairs, and give an example where the pair is not dependent.
Then we define the notion of directionality which deals with counting
the number of coheirs of a type and we give examples of the different
possibilities. Then we discuss non-splintering, an interesting notion
that appears in the work of Rami Grossberg, Andr\'es Villaveces and
Monica VanDieren, and we show that it is not trivial (in the sense
that it can be different than splitting) whenever the directionality
of the theory is not small. In the appendix we study dense types in
RCF. 
\end{abstract}
\maketitle

\section{Introduction}

This paper gives some examples of dependent theories that exemplify
certain phenomenons. Recall,
\begin{defn}
A first order theory $T$ is \emph{dependent} (NIP) if it does not
have the independence property which means: there are no formula $\varphi\left(x,y\right)$
and tuples $\left\langle a_{i},b_{s}\left|\, i<\omega,s\subseteq\omega\right.\right\rangle $
in $\C$ such that $\varphi\left(a_{i},b_{s}\right)$ if and only
if $i\in s$.
\end{defn}

\subsection{Existence of indiscernibles}

Indiscernible sequences are very important in model theory. Usually
one uses Ramsey's theorem to prove their existence. Sometimes we want
to have a stronger result. For instance, we may want that any large
enough set contains an indiscernible sequence and indeed this was
conjectured by Shelah for dependent theories. We will show that at
least in some models of ZFC, one cannot hope for such a result to
be true.

\subsection{Generic pairs}

In a series of papers (\cite{Sh:900,Sh877,Sh906,Sh950}), Shelah has
proved (among other things) that dependent theories give rise to a
``generic pair'' of models (and in fact this characterizes dependent
theories). The natural question is whether the theory of the pair
is again dependent. The answer is no. We present an example of an
$\omega$-stable theory all of whose generic pairs have the independence
property.

\subsection{Directionality}

The directionality of a theory measures the number of finitely satisfiable
global extensions of a complete type (these are also called coheirs).
We say that a theory has \emph{small} directionality if for every
type $p$ over a model $M$, the number of complete finitely satisfiable
(in $M$) $\Delta$-types which are consistent with $p$ is finite
for all finite sets $\Delta$. The theory has \emph{medium} directionality
if this number is bounded by $\left|M\right|$, and it has \emph{large}
directionality if it is not small or medium. We give an equivalent
definition (Theorem \ref{thm:trichotomy} below).

We provide examples of dependent theories of each kind of directionality,
and calculate the directionality of some theories, including RCF and
ACVF.

\subsection{Splintering}

This section is connected to the work of Rami Grossberg, Andr\'es
Villaveces and Monica VanDieren. In \cite{GrViVa} they study Shelah's
Generic pair conjecture (which is now a theorem) and in their analysis,
they came up with the notion of splintering which is similar to splitting.
We show that in any dependent theory with medium or large directionality,
splintering is different than splitting. We also provide an example
of such a theory with small directionality, and prove this cannot
happen in the stable realm.

\subsection{Dense types}

In the appendix, we study dense types in RCF. Namely, we show that
$\ded\lambda$ --- the supremum of the number of cuts of a linear
order of size $\lambda$ --- equals the supremum of the number of
dense types in a model of RCF of size $\lambda$. This is useful for
the calculation of the directionality of RCF.

\subsection{Acknowledgment}

We would like to thank the anonymous referee for many useful comments
and for suggesting to apply the method used for calculating the directionality
of RCF to valued fields (in the previous version it was only shown
that certain valued fields are not small). We would also like to thank
Marcus Tressl with whom we discussed the directionality of RCF and
Pierre Simon and Immanuel Halupczok for discussing valued fields with
us.

\subsection{Notation}

When $\alpha$ and $\beta$ are ordinals, we use left exponentiation
$\leftexp{\beta}{\alpha}$ to denote the set of functions from $\beta$
to $\alpha$, as to not to confuse with ordinal (or cardinal) exponentiation.
If there is no room for confusion, and $A$ and $B$ are some sets
we use $A^{B}$ instead. The set $\alpha^{<\beta}$ is the set of
sequences (functions) $\bigcup\set{\leftexp{\gamma}{\alpha}}{\gamma<\beta}$. 

We do not distinguish elements and tuples unless we say so explicitly.

$\C$ will be the monster model of the theory.

$S_{n}\left(A\right)$ is the set of all complete types in $n$ variables
over $A$. $S_{<\omega}\left(A\right)$ is the union $\bigcup_{n<\omega}S_{n}\left(A\right)$.
$S\left(A\right)$ is the set of all types (perhaps with infinitely
mane variables) over $A$.

For a set of formulas with a partition of variables, $\Delta\left(x,y\right)$,
$L_{\Delta}\left(A\right)$ is the set of formulas of the form $\varphi\left(x,a\right),\neg\varphi\left(x,a\right)$
where $\varphi\left(x,y\right)\in\Delta$ and $a\in A$. $S_{\Delta}\left(A\right)$
is the set of all complete $L_{\Delta}\left(A\right)$-types. Similarly
we may define $\tp_{\Delta}\left(b/A\right)$ as the set of formulas
$\varphi\left(x,a\right)$ such that $\varphi\left(x,y\right)\in\Delta$
and $\C\models\varphi\left(b,a\right)$. For a partial type $p\left(x\right)$
over $A$, $p\upharpoonright\Delta=p\cap L_{\Delta}\left(A\right)$. 

Usually we want to consider a set of formulas $\Delta$ without specifying
a partition of the variables. In this case, for a tuple of variables
$x$, $\Delta^{x}$ is a set of partitioned formulas induced from
$\Delta$ by partitioning the formulas in $\Delta$ to $\left(x,y\right)$
in all possible ways. Then $L_{\Delta}^{x}\left(A\right)$ is just
$L_{\Delta^{x}}\left(A\right)$ and $S_{\Delta}^{x}\left(A\right)$,
$p\upharpoonright\Delta^{x}$ are defined similarly. If $x$ is clear
from the context, we omit it. So for instance, when $p$ is a type
in $x$ over $A$, then $p\upharpoonright\Delta$ is the set of all
formulas $\varphi\left(x,a\right)$ where $\varphi\left(z,w\right)\in\Delta$.

\section{\label{sec:Indisc}Few indiscernibles}

\subsection{Introduction}
\begin{defn}
Let $T$ be a theory. For a cardinal $\kappa$, $n\leq\omega$ and
an ordinal $\delta$, $\kappa\to\left(\delta\right)_{T,n}$ means:
for every set $A\subseteq\C^{n}$ of size $\kappa$, there is a non-constant
sequence of elements of $A$ of length $\delta$ which is indiscernible.
\end{defn}
This definition was suggested by Grossberg and Shelah in \cite[pg. 208, Definition 3.1(2)]{sh:d}
with a slightly different form%
\footnote{The definition there is: $\kappa\to\left(\delta\right)_{T,n}$ if
and only if for each sequence of length $\kappa$ (of $n$-tuples),
there is an indiscernible sub-sequence of length $\delta$. For us
there is no difference because we are dealing with examples where
$\kappa\not\to\left(\mu\right)_{T,n}$. It is also not hard to see
that when $\delta$ is an infinite cardinal these two definitions
are equivalent. %
}. 

There it is also conjectured:
\begin{conjecture}
\label{conj:existence of indiscernibles-1}\cite[pg. 209, Conjecture 3.3]{sh:d}
If $T$ is dependent then for every cardinal $\mu$ there is some
cardinal $\lambda$ such that $\lambda\to\left(\mu\right)_{T,1}$.
\end{conjecture}
In stable theories this holds: it is known that for any $\lambda$
satisfying $\lambda=\lambda^{\left|T\right|}$, $\lambda^{+}\to\left(\lambda^{+}\right)_{T,n}$
(proved by Shelah in \cite{Sh:c}, and follows from local character
of non-forking). In \cite[pg. 209]{sh:d} it is proved that this conjecture
does not hold for simple unstable theories. In \cite{Sh863}, Shelah
proved this conjecture for strongly dependent theories:
\begin{fact}
\label{fac:TStrongly} If $T$ is strongly dependent (see Definition
\ref{def:StronglyDep} below), then for all $\lambda\geq\left|T\right|$,
$\beth_{\left|T\right|^{+}}\left(\lambda\right)\to\left(\lambda^{+}\right)_{T,n}$
for all $n<\omega$.
\end{fact}
This conjuncture is connected to a result by Shelah and Cohen: in
\cite{ShCo919}, they proved that a theory is stable if and only if
it can be presented in some sense in a free algebra\textcolor{black}{{}
in a fixed vocabulary but allowing function symbols with infinite
arity.} If this result could be extended to: a theory is dependent
if and only if it can be represented as an algebra with ordering,
then this could be used to prove existence of indiscernibles.

In this section, we shall show:
\begin{thm}
\label{thm:IntoMain}There is a countable dependent theory $T$ such
that if $\kappa$ is smaller than the first inaccessible cardinal,
then for all $n\in\omega$, $\kappa\not\to\left(\omega\right)_{T,n}$.
\end{thm}
Thus, Conjecture \ref{conj:existence of indiscernibles-1} fails in
a model of ZFC with no inaccessible cardinals. It appears in a more
precise way as Theorem \ref{thm:MainThm} below.

An even stronger result can be obtained, namely:
\begin{fact}
\cite{KaSh975}\label{thm:general result}For every $\theta$ there
is a dependent theory $T$ of size $\theta$ such that for all $\kappa$
and $\delta$, \textup{$\kappa\to\left(\delta\right)_{T,1}$ if and
only if $\kappa\to\left(\delta\right)_{\theta}^{<\omega}$. }
\end{fact}
where:
\begin{defn}
$\kappa\to\left(\delta\right)_{\theta}^{<\omega}$ means: for every
function $c:\left[\kappa\right]^{<\omega}\to\theta$ there is an homogeneous
sub-sequence of length $\delta$ (i.e., there exists $\left\langle \alpha_{i}\left|\, i<\delta\right.\right\rangle \in\leftexp{\delta}{\kappa}$
and $\left\langle c_{n}\left|\, n<\omega\right.\right\rangle \in\leftexp{\omega}{\theta}$
such that $c\left(\alpha_{i_{0}},\ldots,\alpha_{i_{n-1}}\right)=c_{n}$
for every $i_{0}<\cdots<i_{n-1}<\delta$). 
\end{defn}
By \cite{KaSh975}, whenever $\left|T\right|\leq\theta$, $\kappa\to\left(\delta\right)_{\theta}^{<\omega}$
always implies that $\kappa\to\left(\delta\right)_{T,n}$ for all
$n<\omega$, so this is the best result possible. However, the proof
of Theorem \ref{thm:general result} is considerably harder, so it
is given in a subsequent work. 

The second part of this section is devoted to giving a related example
in the field of real numbers. By Fact \ref{fac:TStrongly}, as RCF
is strongly dependent, we cannot prove Theorem \ref{thm:IntoMain}
for RCF, but instead we show that the requirement that $n<\omega$
is necessary:
\begin{thm}
\label{thm:IntoMainRCF}If $\kappa$ is smaller than the first strongly
inaccessible cardinal, then $\kappa\not\to\left(\omega\right)_{RCF,\omega}$.
\end{thm}
This is Theorem \ref{thm:MainThmRCF} below.

\subsubsection*{Notes}

It was unknown to us that in 2011 Kuda{\u\i}bergenov proved a related
result, which refutes a strong version of Conjecture \ref{conj:existence of indiscernibles-1},
namely that $\beth_{\omega+\omega}\left(\mu+\left|T\right|\right)\to\left(\mu\right)_{T,1}$.
He proved that for every ordinal $\alpha$ there exists a dependent
theory (we have not checked whether it is strongly dependent) $T_{\alpha}$
such that $\left|T_{\alpha}\right|=\left|\alpha\right|+\aleph_{0}$
and $\beth_{\alpha}\left(\left|T_{\alpha}\right|\right)\not\to\left(\aleph_{0}\right)_{T_{\alpha},1}$
and thus seem to indicate that the bound in Fact \ref{fac:TStrongly}
is tight. See \cite{russianConjecture}.

\subsubsection*{The idea of the construction}

The counterexample is a ``tree of trees'' with functions connecting
the different trees. For every $\eta$ in the tree $2^{<\omega}$
we shall have a predicate $P_{\eta}$ and an ordering $<_{\eta}$
such that $\left(P_{\eta},<_{\eta}\right)$ is a dense tree. In addition
we shall have functions $G_{\eta,\eta\concat\left\langle i\right\rangle }:P_{\eta}\to P_{\eta\concat\left\langle i\right\rangle }$
for $i=0,1$. The idea is to prove that $\kappa\not\to\left(\mu\right)_{T,1}$
by induction on $\kappa$. To use the induction hypothesis, we push
the counter examples we already have for smaller $\kappa$'s to deeper
levels in the tree $2^{<\omega}$.

\subsection{Preliminaries}
\begin{defn}
We shall need the following fact: \end{defn}
\begin{fact}
\label{fac:DepPolBd} \cite[II, 4]{Sh:c} Let $T$ be any theory.
Then for all $n<\omega$, $T$ is dependent if and only if $\square_{n}$
if and only if $\square_{1}$ where for all $n<\omega$,\end{fact}
\begin{itemize}
\item [$\square_n$] For every finite set of formulas $\Delta\left(x,y\right)$
with $n=\lg\left(x\right)$, there is a polynomial $f$ such that
for every finite set $A\subseteq M\models T$, $\left|S_{\Delta}\left(A\right)\right|\leq f\left(\left|A\right|\right)$.
\end{itemize}
Since we also discuss strongly dependent theories, here is the definition:
\begin{defn}
\label{def:StronglyDep}A theory is called \emph{strongly dependent}
if there is no sequence of formulas $\left\langle \varphi_{n}\left(x,y_{n}\right)\left|\, n<\omega\right.\right\rangle $
such that the set $\left\{ \varphi_{n}\left(x_{\eta},y_{n,k}\right)^{\eta\left(n\right)=k}\left|\,\eta:\omega\to\omega\right.\right\} $
is consistent with the theory (where $\varphi^{\mbox{True}}=\varphi,\varphi^{\mbox{False}}=\neg\varphi$). 
\end{defn}
See \cite{Sh863} for further discussion of strongly dependent theories.
There it is proved that $Th\left(\mathbb{R}\right)$ is strongly dependent,
and so is the theory of the $p$-adics.

\subsection{The example}

Let $S_{n}$ be the finite binary tree $2^{\leq n}$. On a well ordered
tree such as $S_{n}$, we define $<_{\suc}$ as follows: $\eta<_{\suc}\nu$
if $\nu$ is a successor of $\eta$ in the tree.

Let $L_{n}$ be the following language:

\[
L_{n}=\left\{ P_{\eta},<_{\eta},\wedge_{\eta},G_{\eta,\nu}\left|\,\eta,\nu\in S_{n},\eta<_{\suc}\nu\right.\right\} .
\]

Where:
\begin{itemize}
\item $P_{\eta}$ is a unary predicate; $<_{\eta}$ is a binary relation
symbol; $\wedge_{\eta}$ is a binary function symbol; $G_{\eta,\nu}$
is a unary function symbol.
\end{itemize}
Let $T_{n}^{\forall}$ be the following theory:
\begin{itemize}
\item $P_{\eta}\cap P_{\nu}=\emptyset$ for $\eta\neq\nu$. 
\item $\left(P_{\eta},<_{\eta},\wedge_{\eta}\right)$ is a tree, where $\wedge_{\eta}$
is the meet function on $P_{\eta}$, i.e., 
\[
x\wedge_{\eta}y=\max\left\{ z\in P_{\eta}\left|\, z\leq_{\eta}x\,\&\, z\leq_{\eta}y\right.\right\} .
\]

\item $G_{\eta,\nu}:P_{\eta}\to P_{\nu}$ and no further restrictions on
it.
\item In all the axioms above, for elements or pairs outside of the domain
of any of the functions $\wedge_{\eta}$ or $G_{\eta,\nu}$, these
functions are the identity on the leftmost coordinate, so for example
if $\left(x,y\right)\notin P_{\eta}^{2}$, then $x\wedge_{\eta}y=x$. 
\end{itemize}
Thus we have:
\begin{claim}
$T_{n}^{\forall}$ is a universal theory.
\end{claim}

\begin{claim}
$T_{n}^{\forall}$ has the joint embedding property (JEP) and the
amalgamation property (AP). \end{claim}
\begin{proof}
Easy to see. 
\end{proof}
From this we deduce, by e.g., \cite[Theorem 7.4.1]{Hod}:
\begin{cor}
\label{cor:ModelCom}$T_{n}^{\forall}$ has a model completion, $T_{n}$
which eliminates quantifiers, and moreover: if $M\models T_{n+1}^{\forall}$,
$M'=M\upharpoonright L_{n}$ and $M'\subseteq N'\models T_{n}^{\forall}$
then $N'$ can be enriched to a model $N$ of $T_{n+1}^{\forall}$
so that $M\subseteq N$. Hence if $M$ is an existentially closed
model of $T_{n+1}^{\forall}$, then $M'$ is an e.c. model of $T_{n}^{\forall}$.
Hence $T_{n}\subseteq T_{n+1}$ (for more see \cite[Theorem 8.2.4]{Hod}).\end{cor}
\begin{proof}
The moreover part: for each $\eta\in S_{n+1}\backslash S_{n}$, we
define $P_{\eta}^{N}=P_{\eta}^{M}$ and in the same way $\wedge_{\eta}$.
The functions $G_{\eta,\nu}^{N}$ for $\eta\in S_{n}$ and $\nu\in S_{n+1}$
will be extensions of $G_{\eta,\nu}^{M}$. 
\end{proof}
Now we show that $T_{n}$ is dependent, but before that, a few easy
remarks:
\begin{obs}
\label{Obs:Bijection}$ $
\begin{enumerate}
\item If $A\subseteq M\models T_{0}^{\forall}$ is a finite substructure
(so just a tree, with no extra structure), then for all $b\in M$,
the structure generated by $A$ and $b$ is $A\cup\left\{ b\right\} \cup\left\{ \max\left\{ b\wedge a\left|\, a\in A\right.\right\} \right\} $.
\item If $M\models T_{n}^{\forall}$ and $\eta\in2^{\leq n}$, we can define
a new structure $M_{\eta}\models T_{n-\lg\left(\eta\right)}^{\forall}$
whose universe is $\bigcup\left\{ P_{\eta\concat\nu}^{M}\left|\,\nu\in2^{\leq n-\lg\left(n\right)}\right.\right\} $
by: $P_{\nu}^{M_{\eta}}=P_{\eta\concat\nu}^{M}$, and in the same
way we interpret every other symbol (for instance, $G_{\nu_{1},\nu_{2}}^{M_{\eta}}=G_{\eta\concat\nu_{1},\eta\concat\nu_{1}}^{M}$).
For every formula $\varphi\left(x\right)\in L_{n-\lg\left(\eta\right)}$
there is a formula $\varphi'\left(x\right)\in L_{n}$ such that for
all $a\in M_{\eta}$, $M\models\varphi'\left(a\right)$ if and only
if $M_{\eta}\models\varphi\left(a\right)$ (we get $\varphi'$ by
concatenating $\eta$ before any symbol).
\item For $M$ as before and $\eta\in2^{\leq n}$, for any $k<\omega$ there
is a bijection between 
\[
\left\{ p\left(x_{0},\ldots,x_{k-1}\right)\in S_{k}^{\qe}\left(M\right)\left|\,\forall i<k\left(P_{\eta}\left(x_{i}\right)\in p\right)\right.\right\} 
\]
 and 
\[
\left\{ p\left(x_{0},\ldots,x_{k-1}\right)\in S_{k}^{\qe}\left(M_{\eta}\right)\left|\,\forall i<k\left(P_{\left\langle \right\rangle }\left(x_{i}\right)\in p\right)\right.\right\} .
\]

\end{enumerate}
\end{obs}
\begin{proof}
(3): The bijection is given by (2). This is well defined, meaning
that if $p\left(x_{0},\ldots,x_{k-1}\right)$ is a type over $M_{\eta}$
such that $P_{\left\langle \right\rangle }\left(x_{i}\right)\in p$
for all $i<k$, then $\left\{ \varphi'\left|\,\varphi\in p\right.\right\} $
determines a complete type over $M$, such that $P_{\eta}\left(x_{i}\right)\in p$
for all $i<k$. The point is that all atomic formulas over $M$ which
mention elements from $M\backslash M_{\eta}$ or any $\nu\not\geq\eta$
are trivially determined. \end{proof}
\begin{prop}
$T_{n}$ is dependent.\end{prop}
\begin{proof}
We use Fact \ref{fac:DepPolBd}. It is sufficient to find a polynomial
$f\left(x\right)$ such that for every finite set $A$, $\left|S_{1}\left(A\right)\right|\leq f\left(\left|A\right|\right)$. 

First we note that for a set $A$, the size of the structure generated
by $A$ is bounded by a polynomial in $\left|A\right|$: it is generated
by applying $\wedge_{\left\langle \right\rangle }$ on $P_{\left\langle \right\rangle }\cap A$,
applying $G_{\left\langle \right\rangle ,\left\langle 1\right\rangle }$
and $G_{\left\langle \right\rangle ,\left\langle 0\right\rangle }$,
and then applying $\wedge_{\left\langle 0\right\rangle },\wedge_{\left\langle 1\right\rangle }$
and so on. Every step in the process is polynomial, and it ends after
$n$ steps.

Hence we can assume that $A$ is a substructure, i.e., $A\models T_{n}^{\forall}$. 

The proof is by induction on $n$. To ease notation, we shall omit
the subscript $\eta$ from $<_{\eta}$ and $\wedge_{\eta}$.

First we deal with the case $n=0$. In $T_{0}$, $P_{\left\langle \right\rangle }$
is a a tree with no extra structure, while outside $P_{\left\langle \right\rangle }$
there is no structure at all. The number of types outside $P_{\left\langle \right\rangle }$
is bounded by $\left|A\right|+1$ (because there is only one non-algebraic
type). In the case that $P_{\left\langle \right\rangle }\left(x\right)\in p$
for some type $p$ over $A$, we can characterize $p$ by characterizing
the (tree) order-type of $x':=\max\left\{ a\wedge x\left|\, a\in A\right.\right\} $,
i.e., the cut that $x'$ induces on the tree, and by knowing whether
$x'=x$ or $x>x'$ (we note that in general, every theory of a tree
is dependent by \cite{Parigot}). 

Now assume that the claim is true for $n$. Suppose $\eta\in2^{\leq n+1}$
and $1\leq\lg\left(\eta\right)$. By Observation \ref{Obs:Bijection}(3),
there is a bijection between the types $p\left(x\right)$ over $A$
where $P_{\eta}\left(x\right)\in p$ and the types $p\left(x\right)$
in $T_{n+1-\lg\left(\eta\right)}$ over $A_{\eta}$ where $P_{\left\langle \right\rangle }\in p$.
$A_{\eta}\models T_{n+1-\lg\left(\eta\right)}^{\forall}$, and so
by the induction hypothesis, the number of types over $A_{\eta}$
is bounded by a polynomial in $\left|A_{\eta}\right|\leq\left|A\right|$.
As the number of types $p\left(x\right)$ such that $P_{\eta}\left(x\right)\notin p$
for all $\eta$ is bounded by $\left|A\right|+1$ as in the previous
case, we are left with checking the number of types $p\left(x\right)$
such that $P_{\left\langle \right\rangle }\left(x\right)\in p$.

In order to describe $p$, we first have to describe $p$ restricted
to the language $\left\{ <_{\left\langle \right\rangle },\wedge_{\left\langle \right\rangle }\right\} $,
and this is polynomially bounded. Let $x'=\max\left\{ a\wedge x\left|\, a\in A\right.\right\} $.
By Observation \ref{Obs:Bijection}(1), if $A\cup\left\{ x\right\} $
is not closed under $\wedge_{\left\langle \right\rangle }$, $x'$
is the only new element in the structure generated by $A\cup\left\{ x\right\} $
in $P_{\left\langle \right\rangle }$. Hence we are left to determine
the type of the pairs $\left(G_{\left\langle \right\rangle ,\left\langle i\right\rangle }\left(x\right),G_{\left\langle \right\rangle ,\left\langle i\right\rangle }\left(x'\right)\right)$
over $A$ for $i=0,1$ (if $x'$ is not new, then it's enough to determine
the type of $G_{\left\langle \right\rangle ,\left\langle i\right\rangle }\left(x\right)$).
The number of these types is equal to the number of types of pairs
in $T_{n}$ over $A_{\left\langle i\right\rangle }$. As $T_{n}$
is dependent we are done by Fact \ref{fac:DepPolBd}.\end{proof}
\begin{defn}
Let $L=\bigcup_{n<\omega}L_{n}$, $T=\bigcup_{n<\omega}T_{n}$ and
$T^{\forall}=\bigcup_{n<\omega}T_{n}^{\forall}$.
\end{defn}
We easily have:
\begin{cor}
$T$ is complete, it eliminates quantifiers and is dependent.
\end{cor}
We shall prove the following theorem (which implies Theorem \ref{thm:IntoMain}
from the introduction):
\begin{thm}
\label{thm:MainThm}For any two cardinals $\mu\leq\kappa$ such that
in $\left[\mu,\kappa\right]$ there are no (uncountable) strongly
inaccessible cardinals, $\kappa\not\to\left(\mu\right)_{T,1}$.
\end{thm}
We shall prove a slightly stronger statement, by induction on $\kappa$:
\begin{prop}
Given $\mu$ and $\kappa$, such that either $\kappa<\mu$ or there
are no (uncountable) strongly inaccessible cardinals in $\left[\mu,\kappa\right]$,
there is a model $M\models T^{\forall}$ such that $\left|P_{\left\langle \right\rangle }^{M}\right|\geq\kappa$
and $P_{\left\langle \right\rangle }^{M}$ does not contain a non-constant
indiscernible sequence (for quantifier free formulas) of length $\mu$.
\end{prop}
From now on, indiscernible will only mean ``indiscernible for quantifier
free formulas''. 
\begin{proof}
Fix $\mu$. The proof is by induction on $\kappa$. We divide into
cases:
\begin{casenv}
\item $\kappa<\mu$. Clear. 
\item $\kappa=\mu=\aleph_{0}$. Denote $\eta_{j}=\left\langle 1,\ldots,1\right\rangle $,
i.e., the constant sequence of length $j$ and value $1$. Find $M\models T^{\forall}$
such that its universe contains a set $\left\{ a_{i,j}\left|\, i,j<\omega\right.\right\} $
where $a_{i,j}\neq a_{i',j'}$ for all $\left(i,j\right)\neq\left(i',j'\right)$,
$a_{i,j}\in P_{\eta_{j}}$ and in addition $G_{\eta_{j},\eta_{j+1}}\left(a_{i,j}\right)=a_{i,j+1}$
if $j<i$ and $G_{\eta_{j},\eta_{j+1}}\left(a_{i,j}\right)=a_{0,j+1}$
otherwise. We also need that $P_{\left\langle \right\rangle }^{M}=\left\{ a_{i,0}\left|\, i<\omega\right.\right\} $.
Any model satisfying these properties will do (so no need to specify
what the tree structures are). Now, if in $P_{\left\langle \right\rangle }^{M}=\left\{ a_{i,0}\left|\, i<\omega\right.\right\} $
there is a non-constant indiscernible sequence, $\left\langle a_{i_{k},0}\left|\, k<\omega\right.\right\rangle $,
then for $j\geq i_{0},i_{1}$, 
\[
G_{\eta_{j},\eta_{j+1}}\circ\cdots\circ G_{\eta_{0},\eta_{1}}\left(a_{i_{0},0}\right)=G_{\eta_{j},\eta_{j+1}}\circ\cdots\circ G_{\eta_{0},\eta_{1}}\left(a_{i_{1},0}\right).
\]
 But for every $k$ such that $i_{k}>j$, $G_{\eta_{j},\eta_{j+1}}\circ\cdots\circ G_{\eta_{0},\eta_{1}}\left(a_{i_{1},0}\right)\neq G_{\eta_{j},\eta_{j+1}}\circ\cdots\circ G_{\eta_{0},\eta_{1}}\left(a_{i_{k},0}\right)$
--- contradiction.
\item $\kappa$ is singular. Suppose $\kappa=\bigcup_{i<\sigma}\lambda_{i}$
where $\sigma,\lambda_{i}<\kappa$ for all $i<\sigma$. By the induction
hypothesis, for $i<\sigma$ there is a model $M_{i}\models T^{\forall}$
such that $\left|P_{\left\langle \right\rangle }^{M_{i}}\right|\geq\lambda_{i}$
and in $P_{\left\langle \right\rangle }^{M_{i}}$ there is no non-constant
indiscernible sequence of length $\mu$. Also, there is a model $N$
such that $\left|P_{\left\langle \right\rangle }^{N}\right|\geq\sigma$
and in $P_{\left\langle \right\rangle }^{N}$ there is no non-constant
indiscernible sequence of length $\mu$. We may assume that the universes
of all these models are pairwise disjoint and disjoint from $\kappa$.

Suppose that $\left\{ a_{i}\left|\, i<\sigma\right.\right\} \subseteq P_{\left\langle \right\rangle }^{N}$,
and $\left\{ b_{j}\left|\,\sum_{l<i}\lambda_{l}\leq j<\lambda_{i}\right.\right\} \subseteq P_{\left\langle \right\rangle }^{M_{i}}$
witness that $\left|P_{\left\langle \right\rangle }^{N}\right|\geq\sigma$
and $\left|P_{\left\langle \right\rangle }^{M_{i}}\right|\geq\lambda_{i}\backslash\sum_{l<i}\lambda_{l}$.
Let $\bar{M}$ be a model extending each $M_{i}$ and containing the
disjoint union of the sets $\bigcup_{i<\sigma}M_{i}$ (exists by JEP).

Define a new model $M\models T^{\forall}$: $\left(P_{\left\langle \right\rangle }^{M},<_{\left\langle \right\rangle }\right)=\left(\kappa,<\right)$
(so $\wedge_{\left\langle \right\rangle }=\min$); $\left(P_{\left\langle 1\right\rangle \concat\eta}^{M},<_{\eta}\right)=\left(P_{\eta}^{N},<_{\eta}\right)$
and $\left(P_{\left\langle 0\right\rangle \concat\eta}^{M},<_{\left\langle 0\right\rangle ,\eta}\right)=\left(P_{\eta}^{\bar{M}},<_{\eta}\right)$.
In the same way define $\wedge_{\eta}$ for all $\eta$ of length
$\geq1$. The functions are also defined in the same way: $G_{\left\langle 1\right\rangle \concat\eta,\left\langle 1\right\rangle \concat\nu}^{M}=G_{\eta,\nu}^{N}$
and $G_{\left\langle 0\right\rangle \concat\eta,\left\langle 0\right\rangle \concat\nu}^{M}=G_{\eta,\nu}^{\bar{M}}$.
We are left to define $G_{\left\langle \right\rangle ,\left\langle 0\right\rangle }$
and $G_{\left\langle \right\rangle ,\left\langle 1\right\rangle }$.
So let: $G_{\left\langle \right\rangle ,\left\langle 1\right\rangle }\left(\alpha\right)=a_{\min\left\{ i\left|\,\alpha<\lambda_{i}\right.\right\} }$
and $G_{\left\langle \right\rangle ,\left\langle 0\right\rangle }\left(\alpha\right)=b_{\alpha}$
for all $\alpha<\kappa$.

Note that if $I$ is an indiscernible sequence contained in $P_{\left\langle 1\right\rangle }^{M}$
then $I$ is an indiscernible sequence in $N$ contained in $P_{\left\langle \right\rangle }^{N}$,
and the same is true for $P_{\left\langle 0\right\rangle }^{M}$ and
$\bar{M}$.

Assume $\left\langle \alpha_{j}\left|\, j<\mu\right.\right\rangle $
is an indiscernible sequence in $P_{\left\langle \right\rangle }^{M}$.
Then $\left\langle G_{\left\langle \right\rangle ,\left\langle 1\right\rangle }\left(\alpha_{j}\right)\left|\, j<\mu\right.\right\rangle $
is a constant sequence (by the choice of $N$). So there is $i<\sigma$
such that $\sum_{l<i}\lambda_{l}\leq\alpha_{j}<\lambda_{i}$ for all
$j<\mu$. So $\left\langle G_{\left\langle \right\rangle ,\left\langle 0\right\rangle }\left(\alpha_{j}\right)=b_{\alpha_{j}}\left|\, j<\mu\right.\right\rangle $
is a constant sequence (it is indiscernible in $P_{\left\langle \right\rangle }^{\bar{M}}$
and in fact contained in $P_{\left\langle \right\rangle }^{M_{i}}$),
hence $\left\langle \alpha_{j}\left|\, j<\mu\right.\right\rangle $
is constant, as we wanted.

\item $\kappa$ is regular uncountable. By the hypothesis of the proposition,
$\kappa$ is not strongly inaccessible, so there is some $\lambda<\kappa$
such that $2^{\lambda}\geq\kappa$. By the induction hypothesis on
$\lambda$, there is a model $N\models T^{\forall}$ such that in
$P_{\left\langle \right\rangle }^{N}$ there is no non-constant indiscernible
sequence of length $\mu$. Let $\left\{ a_{i}\left|\, i\leq\lambda\right.\right\} \subseteq P_{\left\langle \right\rangle }^{N}$
witness that $\left|P_{\left\langle \right\rangle }^{N}\right|\geq\lambda$.

Define $M\models T^{\forall}$ as follows: $P_{\left\langle \right\rangle }^{M}=2^{\leq\lambda}$
and the ordering is inclusion (equivalently, the ordering is by initial
segment). $\wedge_{\left\langle \right\rangle }$ is defined naturally:
$f\wedge_{\left\langle \right\rangle }g=f\upharpoonright\min\left\{ \alpha\left|\, f\left(\alpha\right)\neq g\left(\alpha\right)\right.\right\} $.

For all $\eta$, let $P_{\left\langle 1\right\rangle \concat\eta}^{M}=P_{\eta}^{N}$,
and the ordering and the functions are naturally induced from $N$.
The main point is that we set $G_{\left\langle \right\rangle ,\left\langle 1\right\rangle }\left(f\right)=a_{\lg\left(f\right)}$.
Now choose $P_{\left\langle 0\right\rangle \concat\eta}^{M}$, $G_{\left\langle 0\right\rangle \concat\eta,\left\langle 0\right\rangle \concat\nu}$,
etc. arbitrarily, and let $G_{\left\langle \right\rangle ,\left\langle 0\right\rangle }$
be any function.

Suppose that $\left\langle f_{i}\left|\, i<\mu\right.\right\rangle $
is a non-constant indiscernible sequence:

If $f_{1}<f_{0}$ (i.e., $f_{1}<_{\left\langle \right\rangle }f_{0}$),
we shall have an infinite decreasing sequence in a well-ordered tree
--- a contradiction.

If $f_{0}<f_{1}$, $\left\langle f_{i}\left|\, i<\mu\right.\right\rangle $
is increasing, so $\left\langle G_{\left\langle \right\rangle ,\left\langle 1\right\rangle }^{M}\left(f_{i}\right)=a_{\lg\left(f_{i}\right)}\left|\, i<\mu\right.\right\rangle $
is non-constant --- contradiction (as it is an indiscernible sequence
in $M$ and hence in $P_{\left\langle 0\right\rangle }^{N}$).

Let $h_{i}=f_{0}\wedge f_{i+1}$ for $i<\mu$ (where $\wedge=\wedge_{\left\langle \right\rangle }$).
This is an indiscernible sequence, and by the same arguments, it cannot
increase or decrease, but as $h_{i}\leq f_{0}$, and $\left(P_{\left\langle \right\rangle },<_{\left\langle \right\rangle }\right)$
is a tree, it follows that $h_{i}$ is constant.

Assume $f_{0}\wedge f_{1}<f_{1}\wedge f_{2}$, then $f_{2i}\wedge f_{2i+1}<f_{2\left(i+1\right)}\wedge f_{2\left(i+1\right)+1}$
for all $i<\mu$, and again $\left\langle f_{2i}\wedge f_{2i+1}\left|\, i<\mu\right.\right\rangle $
an increasing indiscernible sequence and we have a contradiction.

By the same reasoning, it cannot be that $f_{0}\wedge f_{1}>f_{1}\wedge f_{2}$.
As $\left(P_{\left\langle \right\rangle },<_{\left\langle \right\rangle }\right)$
is a tree, we conclude that $f_{0}\wedge f_{2}=f_{0}\wedge f_{1}=f_{1}\wedge f_{2}$.
But that is a contradiction (because if $\alpha=\lg\left(f_{0}\wedge f_{1}\right)$,
then $\left|\left\{ f_{0}\left(\alpha\right),f_{1}\left(\alpha\right),f_{2}\left(\alpha\right)\right\} \right|=3$). 

\end{casenv}
\end{proof}

\subsection{In RCF there are few indiscernibles of $\omega$-tuples.}

Here we will prove Theorem \ref{thm:IntoMainRCF}. Since RCF is strongly
dependent, Fact \ref{fac:TStrongly} (which discusses finite tuples)
holds for it, so we will show that a similar phenomenon as in the
previous section holds for $\omega$-tuples in RCF. So assume $\C\models RCF$.
\begin{notation}
The set of all open intervals $\left(a,b\right)$ (where $a<b$ and
$a,b\in\C$) is denoted by $\I$. \end{notation}
\begin{defn}
\label{def:arrowInterval}For a cardinal $\kappa$, $n\leq\omega$
and an ordinal $\delta$, $\kappa\to\left(\delta\right)_{n}^{\intr}$
means: for every set $A$ of $n$-tuples of (non-empty, open) intervals
(so for each $\bar{I}\in A$, $\bar{I}=\left\langle I^{i}\left|\, i<n\right.\right\rangle \in\mathfrak{J}^{n}$)
of size $\kappa$, there is a sequence $\left\langle \bar{I}_{\alpha}\left|\,\alpha<\delta\right.\right\rangle \in A^{\delta}$
of order type $\delta$ such that $\bar{I}_{\alpha}\neq\bar{I}_{\beta}$
for $\alpha<\beta<\delta$, and there is a sequence $\left\langle \bar{b}_{\alpha}\left|\,\alpha<\delta\right.\right\rangle $
such that $\bar{b}_{\alpha}\in\bar{I}_{\alpha}$ (i.e., $\bar{b}_{\alpha}=\left\langle b_{\alpha}^{0},\ldots,b_{\alpha}^{n-1}\right\rangle $
and $b_{\alpha}^{i}\in I_{\alpha}^{i}$) and such that $\left\langle \bar{b}_{\alpha}\left|\,\alpha<\delta\right.\right\rangle $
is an indiscernible sequence. \end{defn}
\begin{rem}
Note that: 
\begin{enumerate}
\item If $\kappa\to\left(\delta\right)_{n}^{\intr}$ then $\kappa\to\left(\delta\right)_{m}^{\intr}$
for all $m\leq n$.
\item If $\kappa\not\to\left(\delta\right)_{n}^{\intr}$ then $\kappa\not\to\left(\delta\right)_{RCF,n}$
(why? if $A$ witnesses that $\kappa\not\to\left(\delta\right)_{n}^{\intr}$,
then for each $\bar{I}\in A$, choose $\bar{b}_{\bar{I}}\in\bar{I}$
(as above) in such a way that $\left\{ \bar{b}_{\bar{I}}\left|\,\bar{I}\in A\right.\right\} $
has size $\kappa$. By definition this set witnesses $\kappa\not\to\left(\delta\right)_{RCF,n}$). 
\item If $\lambda<\kappa$ and $\kappa\not\to\left(\delta\right)_{n}^{\intr}$
then $\lambda\not\to\left(\delta\right)_{n}^{\intr}$. 
\end{enumerate}
\end{rem}
We shall prove the following theorem (which immediately implies Theorem
\ref{thm:IntoMainRCF}):
\begin{thm}
\label{thm:MainThmRCF}For any two cardinals $\mu\leq\kappa$ such
that in $\left[\mu,\kappa\right]$ there are no strongly inaccessible
cardinals, $\kappa\not\to\left(\mu\right)_{\omega}^{\intr}$.
\end{thm}
The proof follows from a sequence of claims:
\begin{claim}
\label{cla:Obvious}If $\kappa<\mu$ then $\kappa\not\to\left(\mu\right)_{n}^{\intr}$
for all $n\leq\omega$.\end{claim}
\begin{proof}
Obvious.\end{proof}
\begin{claim}
\label{cla:aleph0}If $\kappa=\mu=\aleph_{0}$ then $\kappa\not\to\left(\mu\right)_{1}^{\intr}$.\end{claim}
\begin{proof}
For $n<\omega$, let $I_{n}=\left(n,n+1\right)$. \end{proof}
\begin{claim}
\label{cla:Singular}Suppose $\kappa=\sum_{i<\sigma}\lambda_{i}$
and $n\leq\omega$. Then, if $\sigma\not\to\left(\mu\right)_{n}^{\intr}$
and $\lambda_{i}\not\to\left(\mu\right)_{n}^{\intr}$ then $\kappa\not\to\left(\mu\right)_{2+2n}^{\intr}$.\end{claim}
\begin{proof}
By assumption, we have a set of intervals $\left\{ \bar{R}_{i}\left|\, i<\sigma\right.\right\} $
that witness $\sigma\not\to\left(\mu\right)_{n}^{\intr}$ and for
each $i<\sigma$ we have $\left\{ \bar{S}_{\beta}\left|\,\sum_{j<i}\lambda_{j}<\beta<\lambda_{i}\right.\right\} $
that witness $\left|\lambda_{i}\backslash\sum_{j<i}\lambda_{j}\right|\not\to\left(\mu\right)_{n}^{\intr}$. 

Fix an increasing sequence of elements $\left\langle b_{i}\left|\, i<\sigma\right.\right\rangle $.

For $\alpha<\kappa$, let $\beta=\beta\left(\alpha\right)=\min\left\{ i<\sigma\left|\,\alpha<\lambda_{i}\right.\right\} $
and for $i<2+2n$, define:
\begin{itemize}
\item If $i=0$, let $I_{\alpha}^{i}=\left(b_{2\beta\left(\alpha\right)},b_{2\beta\left(\alpha\right)+1}\right)$.
\item If $i=1$, let $I_{\alpha}^{i}=\left(b_{2\beta\left(\alpha\right)+1},b_{2\beta\left(\alpha\right)+2}\right)$. 
\item If $i=2k+2$, let $I_{\alpha}^{i}=R_{\beta\left(\alpha\right)}^{k}$.
\item If $i=2k+3$, let $I_{\alpha}^{i}=S_{\alpha}^{k}$. 
\end{itemize}
Suppose $\left\langle \bar{b}_{\varepsilon}\left|\,\varepsilon<\mu\right.\right\rangle $
is an indiscernible sequence such that $\bar{b}_{\varepsilon}\in\bar{I}_{\alpha_{\varepsilon}}$
for $\varepsilon<\mu$. Denote $\bar{b}_{\varepsilon}=\left\langle b_{0}^{\varepsilon},\ldots,b_{2+2n-1}^{\varepsilon}\right\rangle $.
Note that $b_{1}^{\varepsilon}<b_{0}^{\varepsilon'}$ if and only
if $\beta\left(\alpha_{\varepsilon}\right)<\beta\left(\alpha_{\varepsilon'}\right)$
(we need two intervals for the ``only if'' direction).

Hence $\left\langle \beta\left(\alpha_{\varepsilon}\right)\left|\,\varepsilon<\mu\right.\right\rangle $
is increasing or constant. But if it is increasing then we have a
contradiction to the choice of $\left\{ \bar{R}_{i}\left|\, i<\sigma\right.\right\} $.
So it is constant, and suppose $\beta\left(\alpha_{\varepsilon}\right)=i_{0}$
for all $\varepsilon<\mu$. But then $\alpha_{\varepsilon}\in\lambda_{i_{0}}\backslash\sum_{j<i_{0}}\lambda_{j}$
for all $\varepsilon<\mu$ and we get a contradiction to the choice
of $\left\langle \bar{S}_{\beta}\left|\,\sum_{j<i_{0}}\lambda_{j}<\beta<\lambda_{i_{0}}\right.\right\rangle $. \end{proof}
\begin{claim}
\label{cla:Regular}Suppose $\lambda\not\to\left(\mu\right)_{n}^{\intr}$.
Then $2^{\lambda}\not\to\left(\mu\right)_{4+2n}^{\intr}$. \end{claim}
\begin{proof}
\renewcommand{\qedsymbol}{}Suppose $\left\{ \bar{I}_{\alpha}\left|\,\alpha<\lambda\right.\right\} $
witnesses that $\lambda\not\to\left(\mu\right)_{n}^{\intr}$. 

By adding two intervals to each $\bar{I}_{\alpha}$, we can ensure
that it has the extra property that if $\bar{c}_{1}\in I_{\alpha_{1}}$
and $\bar{c}_{2}\in\bar{I}_{\alpha_{2}}$ then $c_{1}^{1}<c_{2}^{0}$
if and only if $\alpha_{1}<\alpha_{2}$ (as in the previous claim).
By this we have increased the length of $\bar{I}_{\alpha}$ to $2+n$
(and it is still a witness of $\lambda\not\to\left(\mu\right)_{n}^{\intr}$).

We write $\bar{c}_{1}<^{*}\bar{c}_{2}$ for $c_{1}^{1}<c_{2}^{0}$
, but note that it is not really an ordering (it is not transitive
in general).

We shall find below a four-place definable function $f$ such that: 
\begin{itemize}
\item [$\heartsuit$] For every two ordinals, $\delta,\zeta$, if $\left\langle \bar{R}_{\alpha}\left|\,\alpha<\delta\right.\right\rangle $
is a sequence of $\zeta$-tuples of intervals, then there exists a
set of $2\zeta$-tuples of intervals, $\left\{ \bar{S}_{\eta}\left|\,\eta\in\leftexp{\delta}2\right.\right\} $
(of size $2^{\left|\delta\right|}$) such that for all $i<\zeta$
and $\eta_{1}\neq\eta_{2}$, if $b_{1}\in S_{\eta_{1}}^{2i},b_{2}\in S_{\eta_{1}}^{2i+1}$
and $b_{3}\in S_{\eta_{2}}^{2i},b_{4}\in S_{\eta_{2}}^{2i+1}$ then
$f\left(b_{1},b_{2},b_{3},b_{4}\right)$ is in $R_{\lg\left(\eta_{1}\wedge\eta_{2}\right)}^{i}$. 
\end{itemize}
Apply $\heartsuit$ to our situation to get $\left\{ \bar{J}_{\eta}\left|\,\eta\in\leftexp{\lambda}2\right.\right\} $
such that $\bar{J}_{\eta}=\left\langle J_{\eta}^{i}\left|\, i<4+2n\right.\right\rangle $
and for all $k<2+n$ and $\eta_{1}\neq\eta_{2}$, if $b_{1}\in J_{\eta_{1}}^{2k},b_{2}\in J_{\eta_{1}}^{2k+1}$
and $b_{3}\in J_{\eta_{2}}^{2k},b_{4}\in J_{\eta_{2}}^{2k+1}$ then
$f\left(b_{1},b_{2},b_{3},b_{4}\right)$ is in $I_{\lg\left(\eta_{1}\wedge\eta_{2}\right)}^{k}$. 

This is enough (the reasons are exactly as in the regular case of
the proof of Theorem \ref{thm:MainThm}, but we shall repeat it for
clarity):

To simplify notation, we regard $f$ as a function on tuples, so that
if $\bar{b}_{1}\in\bar{J}_{\eta_{1}},\bar{b}_{2}\in\bar{J}_{\eta_{2}}$
then $f\left(\bar{b}_{1},\bar{b}_{2}\right)$ is in $\bar{I}_{\lg\left(\eta_{1}\wedge\eta_{2}\right)}$
(namely, $f\left(\bar{b}_{1},\bar{b}_{2}\right)=\left\langle a_{k}\left|\, k<2+n\right.\right\rangle $
where $a_{k}=f\left(b_{2k}^{1},b_{2k+1}^{1},b_{2k}^{2},b_{2k+1}^{2}\right)\in I_{\lg\left(\eta_{1}\wedge\eta_{2}\right)}^{k}$
for $k<2+n$).

Suppose $\left\langle \eta_{i}\left|\, i<\mu\right.\right\rangle \subseteq\leftexp{\lambda}2$
is without repetitions and $\left\langle \bar{b}_{\eta_{i}}\left|\, i<\mu\right.\right\rangle $
is an indiscernible sequence such that $\bar{b}_{\eta_{i}}\in\bar{J}_{\eta_{i}}$.

Let $h_{i}=\eta_{0}\wedge\eta_{i+1}$ for $i<\mu$. If $\lg\left(h_{i}\right)<\lg\left(h_{j}\right)$
for some $i\neq j$ then $f\left(\bar{b}_{\eta_{0}},\bar{b}_{\eta_{i+1}}\right)<^{*}f\left(\bar{b}_{\eta_{0}},\bar{b}_{\eta_{j+1}}\right)$
and so by indiscernibility, $\left\langle \lg\left(h_{i}\right)\left|\, i<\mu\right.\right\rangle $
is increasing (it cannot be decreasing), and so $f\left(\bar{b}_{\eta_{0}},\bar{b}_{\eta_{i+1}}\right)$
contradicts our choice of $\left\langle \bar{I}_{\alpha}\left|\,\alpha<\lambda\right.\right\rangle $.
Hence (because $h_{i}\leq\eta_{0}$) $h_{i}$ is constant.

Assume $\eta_{0}\wedge\eta_{1}<\eta_{1}\wedge\eta_{2}$, then $f\left(\bar{b}_{\eta_{0}},\bar{b}_{\eta_{1}}\right)<^{*}f\left(\bar{b}_{\eta_{1}},\bar{b}_{\eta_{2}}\right)$
so $f\left(\bar{b}_{\eta_{1}},\bar{b}_{\eta_{2}}\right)<^{*}f\left(\bar{b}_{\eta_{2}},\bar{b}_{\eta_{3}}\right)$,
and so $\lg\left(\eta_{0}\wedge\eta_{1}\right)<\lg\left(\eta_{1}\wedge\eta_{2}\right)<\lg\left(\eta_{2}\wedge\eta_{3}\right)$
hence, $f\left(\bar{b}_{\eta_{0}},\bar{b}_{\eta_{1}}\right)<^{*}f\left(\bar{b}_{\eta_{2}},\bar{b}_{\eta_{3}}\right)$
and it follows that $\left\langle \lg\left(\eta_{2i}\wedge\eta_{2i+1}\right)\left|\, i<\mu\right.\right\rangle $
is increasing. And this is again a contradiction. 

Similarly, it cannot be that $\eta_{0}\wedge\eta_{1}>\eta_{1}\wedge\eta_{2}$.
As both sides are less or equal than $\eta_{1}$, it must be that
$\eta_{0}\wedge\eta_{2}=\eta_{0}\wedge\eta_{1}=\eta_{1}\wedge\eta_{2}$.
But that is impossible (because if $\alpha=\lg\left(\eta_{0}\wedge\eta_{1}\right)$,
then $\left|\left\{ \eta_{0}\left(\alpha\right),\eta_{1}\left(\alpha\right),\eta_{2}\left(\alpha\right)\right\} \right|=3$).
\begin{claim*}
$\heartsuit$ is true. \end{claim*}
\begin{proof}
\renewcommand{\qedsymbol}{$\square$}Let $f\left(x,y,z,w\right)=\left(x-z\right)/\left(y-w\right)$
(do not worry about division by $0$, we shall explain below).

It is enough, by the definition of $\heartsuit$, to assume $\zeta=1$.
By compactness, we may assume that $\delta$ is finite, and to avoid
confusion, denote it by $m$. So we have a finite set, $\leftexp m2$,
and a sequence of intervals $\left\langle R_{i}\left|\, i<m\right.\right\rangle $.
Each $R_{i}$ is of the form $\left(a_{i},b_{i}\right)$. Let $c_{i}=\left(b_{i}+a_{i}\right)/2$.
Let $d\in\C$ be any element greater than any member of $A:=\acl\left(a_{i},b_{i}\left|\, i<m\right.\right)$.
For each $\eta\in\leftexp m2$, let $a_{\eta}=\sum_{i<m}\eta\left(i\right)c_{i}d^{m-i}$,
and $b_{\eta}=\sum_{i<m}\eta\left(i\right)d^{m-i}$. 

Let $S_{\eta}^{0}=\left(a_{\eta}-1,a_{\eta}+1\right)$ and $S_{\eta}^{1}=\left(b_{\eta}-1,b_{\eta}+1\right)$. 

This works:

Assume that $\eta_{1}\neq\eta_{2}$, $b_{1}\in S_{\eta_{1}}^{0},b_{2}\in S_{\eta_{1}}^{1}$
and $b_{3}\in S_{\eta_{2}}^{0},b_{4}\in S_{\eta_{2}}^{1}$. 

We have to show $\left(b_{1}-b_{3}\right)/\left(b_{2}-b_{4}\right)\in R_{\lg\left(\eta_{1}\wedge\eta_{2}\right)}$.
Denote $k=\lg\left(\eta_{1}\wedge\eta_{2}\right)$ (so $k<m$).

$a_{\eta_{1}}-a_{\eta_{2}}$ is of the form $\varepsilon c_{k}d^{m-k}+F\left(d\right)$
where $\varepsilon\in\left\{ -1,1\right\} $, and $F\left(d\right)$
is a polynomial over $A$ of degree $\leq m-k-1$. $b_{\eta_{1}}-b_{\eta_{2}}$
is of the form $\varepsilon d^{m-k}+G\left(d\right)$, where $\varepsilon$
is the same for both (and $G$ is a polynomial over $\mathbb{Z}$
of degree $\leq m-k-1$). Now, $b_{1}-b_{3}\in\left(a_{\eta_{1}}-a_{\eta_{2}}-2,a_{\eta_{1}}-a_{\eta_{2}}+2\right)$,
and $b_{2}-b_{4}\in\left(b_{\eta_{1}}-b_{\eta_{2}}-2,b_{\eta_{1}}-b_{\eta_{2}}+2\right)$,
and hence we know that $b_{2}-b_{4}\neq0$. It follows that $\left(b_{1}-b_{3}\right)/\left(b_{2}-b_{4}\right)$
is inside an interval whose endpoints are $\left\{ \left(\varepsilon c_{k}d^{m-k}+F\left(d\right)\pm2\right)/\left(\varepsilon d^{m-k}+G\left(d\right)\pm2\right)\right\} $.
But 
\[
\left(\varepsilon c_{k}d^{m-k}+F\left(d\right)\pm2\right)/\left(\varepsilon d^{m-k}+G\left(d\right)\pm2\right)\in R_{k}
\]
by our choice of $d$, and we are done.

Note that for $\eta_{1}\neq\eta_{2}$, $S_{\eta_{1}}^{0}\neq S_{\eta_{2}}^{0}$
regardless of the $R_{i}$'s (which can be a constant interval). 
\end{proof}
\end{proof}
The proof of Theorem \ref{thm:MainThmRCF} now follows by induction
on $\kappa$: fix $\mu$, and let $\kappa$ be the first cardinal
for which the theorem fails. Then by Claim \ref{cla:Obvious}, $\kappa\geq\mu$.
By Claim \ref{cla:aleph0}, $\aleph_{0}<\kappa$. By Claim \ref{cla:Singular},
$\kappa$ cannot be singular. By Claim \ref{cla:Regular}, $\kappa$
cannot be regular, because if it were, there would be a $\lambda<\kappa$
such that $2^{\lambda}\geq\kappa$ (because $\kappa$ is not strongly
inaccessible). Note that we did use Claim \ref{cla:Obvious} to deal
with cases where we couldn't use the induction hypothesis (for example,
in the regular case, it might be that $\lambda<\mu$).

\subsubsection*{Further remarks}

Theorem \ref{thm:MainThmRCF} can be generalized to allow parameters: 

Suppose $\C\models RCF$, and $A\subseteq\C$. 
\begin{defn}
$\kappa\to_{A}\left(\mu\right)_{\omega}^{\intr}$ means the same as
in Definition \ref{def:arrowInterval}, but we require that the indiscernible
sequence is indiscernible over $A$.
\end{defn}
Then we have:
\begin{thm}
For any set of parameters $A$ and any two cardinals $\mu,\kappa$
such that in $\left[\max\left\{ \left|A\right|,\mu\right\} ,\kappa\right]$
there are no strongly inaccessible cardinals or $\kappa<\max\left\{ \left|A\right|,\mu\right\} $,
$\kappa\not\to_{A}\left(\mu\right)_{\omega}^{\intr}$.\end{thm}
\begin{proof}
The proof goes exactly as the proof of Theorem \ref{thm:MainThmRCF},
but the base case for the induction is different. If $\max\left\{ \left|A\right|,\mu\right\} =\mu$,
the proof is exactly the same. Otherwise, we have to deal with the
case $\kappa\leq\left|A\right|$: 

Enumerate $A=\left\{ a_{i}\left|\, i<\mu'\right.\right\} $. Let $\varepsilon\in\C$
be greater than $0$ but smaller than any element in $\acl\left(A\right)$.
For $i<\mu'$, let $I_{i}=\left(a_{i},a_{i}+\varepsilon\right)$.
Then $\left\{ I_{i}\left|\, i<\kappa\right.\right\} $ witnesses $\kappa\not\to_{A}\left(\mbox{\ensuremath{\mu}}\right)_{\omega}^{\intr}$.
\end{proof}

\section{Generic pair}

Here we give an example of an $\omega$-stable theory, such that for
all weakly generic pairs of structures $M\prec M_{1}$ (see below
for the definition) the theory of the pair $\left(M_{1},M\right)$
in an extended language where we name $M$ by a predicate has the
independence property.
\begin{defn}
\label{def:WeaklyGenPair}A pair $\left(M_{1},M\right)$ as above
is \emph{weakly generic} if for all formula $\varphi\left(x\right)$
with parameters from $M$, if $\varphi$ has infinitely many solutions
in $M$, then it has a solution in $M_{1}\backslash M$.
\end{defn}
This definition is induced by the well known ``generic pair conjecture''
(see \cite{Sh:900,Sh950}), and it is worth while to give the precise
definitions.
\begin{defn}
\label{def:GenPair}Assume that $\lambda=\lambda^{<\lambda}>\left|T\right|$
(in particular, $\lambda$ is regular) and that $2^{\lambda}=\lambda^{+}$.
The \emph{generic pair property }for $\lambda$ says that there exists
a saturated model $M$ of cardinality $\lambda^{+}$, an increasing
continuous sequence of models $\sequence{M_{\alpha}}{\alpha<\lambda^{+}}$
and a club $E\subseteq\lambda^{+}$ such that $\bigcup_{\alpha<\lambda^{+}}M_{\alpha}=M$
and for all $\alpha<\beta\in E$ of cofinality $\lambda$, the pair
$\left(M_{\beta},M_{\alpha}\right)$ has the same isomorphism type.
We call this pair the \emph{generic} \emph{pair} of $T$ of size $\lambda$.\end{defn}
\begin{prop}
\label{prop:generic pair depenends only on T}Assume that $\lambda=\lambda^{<\lambda}>\left|T\right|$
and that $2^{\lambda}=\lambda^{+}$. The generic pair property for
$\lambda$ holds iff for every saturated model $N$ of cardinality
$\lambda^{+}$ and for every increasing continuous sequence of models
$\sequence{N_{\alpha}}{\alpha<\lambda^{+}}$ with union $N$ there
exists a club $E\subseteq\lambda^{+}$ such that for all $\alpha<\beta\in E$
of cofinality $\lambda$, the pair $\left(N_{\beta},N_{\alpha}\right)$
has the same isomorphism type. Moreover, this type does not depended
on the particular choice of $N$ or $N_{\alpha}$.\end{prop}
\begin{proof}
Left to right:

Suppose $M$, $\sequence{M_{\alpha}}{\alpha<\lambda^{+}}$ and $E$
witness that $T$ has the generic pair property for $\lambda$. If
$N$ is another saturated model of size $\lambda^{+}$ and $\sequence{N_{\alpha}}{\alpha<\lambda^{+}}$
is as in the Proposition. Then $N\cong M$, so we may assume $N=M$.
Let $E_{0}=\left\{ \delta<\lambda^{+}\left|\, N_{\delta}=M_{\delta}\right.\right\} $.
This is a club of $\lambda^{+}$, and so $E\cap E_{0}$ is also a
club of $\lambda^{+}$ such that $\left(N_{\beta},N_{\alpha}\right)$
has the same isomorphism type for any $\alpha<\beta\in E\cap E_{0}$
of cofinality $\lambda$. 

Right to left is clear. 
\end{proof}
Justifying definition \ref{def:WeaklyGenPair} we have:
\begin{claim}
Assume that $T$ has the generic pair property for $\lambda$, then
every generic pair of size $\lambda$ is weakly generic.\end{claim}
\begin{proof}
Suppose that $M$, $\left\langle M_{\alpha}\left|\,\alpha<\lambda^{+}\right.\right\rangle $
and $E$ are as in Definition \ref{def:GenPair}. Suppose $\alpha,\beta\in E$
and $\alpha<\beta$ are of cofinality $\lambda$. We are given a formula
$\varphi\left(x\right)$ with parameter from $M_{\alpha}$, such that
$\aleph_{0}\leq\left|\varphi\left(M_{\alpha}\right)\right|$. By saturation
of $M$, $\lambda^{+}=\left|\varphi\left(M\right)\right|$. Since
$M=\bigcup_{\beta'\in E,\mbox{cf}\left(\beta'\right)=\lambda}M_{\beta}$
there is some $\alpha<\beta'\in E$ of cofinality $\lambda$ such
that $\varphi\left(M_{\beta'}\right)\backslash M_{\alpha}\neq\emptyset$,
but as $\left(M_{\beta},M_{\alpha}\right)\cong\left(M_{\beta'},M_{\alpha}\right)$,
we are done.
\end{proof}
Proposition \ref{prop:generic pair depenends only on T} implies that
the generic pair property and the the generic pair are both natural
notions. It is important in the study of dependent theories as it
lead to the development of a theory of type decomposition in NIP.
Using this theory, the second author's \cite{Sh950,Sh:900} prove
that the generic pair property holds for dependent theories and large
enough $\lambda$'s. On the other hand, \cite{Sh877,Sh906} prove
that if $T$ has IP then it lacks the generic pair property for all
large enough $\lambda$. 

Hence it makes sense to ask whether the theory of the pair is dependent. 

The answer is no:
\begin{thm}
There exists an $\omega$-stable theory such that for every weakly
generic pair of models $M\prec M_{1}$, the theory of the pair $\left(M_{1},M\right)$
has the independence property.
\end{thm}
We shall describe this theory:

Let $L=\left\{ P,R,Q_{1},Q_{2}\right\} $ where $R,P$ are unary predicates
and $Q_{1},Q_{2}$ are binary relations. 

Let $\tilde{M}$ be the following structure for $L$: 
\begin{enumerate}
\item The universe is: 
\begin{eqnarray*}
\tilde{M} & = & \left\{ u\subseteq\omega\left|\,\left|u\right|<\omega\right.\right\} \cup\\
 &  & \left\{ \left(u,v,i\right)\left|\, u,v\subseteq\omega,\left|u\right|<\omega,\left|v\right|<\omega,i<\omega\,\&\, u\subseteq v\Rightarrow i<\left|v\right|+1\right.\right\} .
\end{eqnarray*}

\item The predicates are interpreted as follows:

\begin{itemize}
\item $P^{\tilde{M}}=\left\{ u\subseteq\omega\left|\,\left|u\right|<\aleph_{0}\right.\right\} $. 
\item $R^{\tilde{M}}$ is $\tilde{M}\backslash\left(P^{M}\right)$.
\item $Q_{1}^{\tilde{M}}=\left\{ \left(u,\left(u,v,i\right)\right)\left|\, u\in P^{\tilde{M}}\right.\right\} $.
\item $Q_{2}^{\tilde{M}}=\left\{ \left(v,\left(u,v,i\right)\right)\left|\, v\in P^{\tilde{M}}\right.\right\} $.
\end{itemize}
\end{enumerate}
Let $T=Th\left(\tilde{M}\right)$. 

As we shall see in the next claim, $T$ gives rise to the following
definition:
\begin{defn}
\label{def:PBA}We call a structure $\left(B,\cup,\cap,-,\subseteq,0\right)$
a \emph{pseudo Boolean algebra} (PBA) when it satisfies all the axioms
of a Boolean algebra except: 

There is no greatest element $1$ (i.e., remove all the axioms concerning
it).
\end{defn}
Pseudo Boolean algebra can have atoms like in Boolean algebras (nonzero
elements that do not contain any smaller nonzero elements). 
\begin{defn}
Say that a PBA is of \emph{finite type} if every element is a union
of finitely many atoms. 
\end{defn}

\begin{defn}
For a PBA $A$, and $C\subseteq A$ a sub-PBA, let $A_{C}:=\left\{ a\in A\left|\,\exists c\in C\left(a\subseteq^{A}c\right)\right.\right\} $,
and for a subset $D\subseteq A$, let $\atom\left(D\right)$ be the
set of atoms contained in $D$.\end{defn}
\begin{prop}
\label{pro:ISomPBA}Every PBA of finite type is isomorphic to $\left(\SS_{<\infty}\left(\kappa\right),\cup,\cap,-,\subseteq,\emptyset\right)$
for some $\kappa$ where $\SS_{<\infty}\left(\kappa\right)$ is the
set of all finite subsets of $\kappa$. Moreover: Assume $A,B$ are
PBAs of finite type and $C\subseteq A,B$ is a common sub-PBA. Then,
if: 
\begin{enumerate}
\item $\left|\atom\left(A\right)\backslash\atom\left(A_{C}\right)\right|=\left|\atom\left(B\right)\backslash\atom\left(B_{C}\right)\right|$.
\item For every $c\in C$, $A$ and $B$ agree on the size of $c$ (the
number of atoms it contains).
\end{enumerate}
Then there is an isomorphism of PBAs $f:A\to B$ such that $f\upharpoonright C=\id$.\end{prop}
\begin{proof}
The first part follows from the easy observation that in a PBA of
finite type, every element has a unique presentation as a union of
finitely many atoms. So if $A$ is a PBA, and its set of atoms is
$\left\{ a_{i}\left|\, i<\kappa\right.\right\} $, then take $a_{i}$
to $\left\{ i\right\} $.

For the moreover part, first we extend $\id_{C}$ to an isomorphism
from $A_{C}$ to $B_{C}$: consider all elements in $C$ of minimal
size, these are the atoms of $C$. For each such $c$, map the set
of atoms in $A$ contained in $c$ to the set of atoms in $B$ contained
in $c$. This is well defined and can be extended to all of $A_{C}$. 

Now, $\left|\atom\left(A\right)\backslash\atom\left(A_{C}\right)\right|=\left|\atom\left(B\right)\backslash\atom\left(B_{C}\right)\right|$,
so any bijection between the set of atoms induces an isomorphism. \end{proof}
\begin{claim}
\label{cla:omegaStable}$T$ is $\omega$-stable. \end{claim}
\begin{proof}
We prove that an expansion of $T$ to a larger vocabulary is $\omega$-stable,
by adding new relations to the language, which are all definable ---
\[
\left\{ S_{n},\subseteq_{n},\pi_{1},\pi_{2},\cap_{n},\cup_{n},-_{n},e\left|\, n\geq1\right.\right\} 
\]
where $S_{n}$ is a unary relation defined on $P$, $\subseteq_{n}$
is a binary relation defined on $P$, $\pi_{1},\pi_{2}$ are two unary
functions from $R$ to $P$, $\cap_{n},-_{n}$ are binary functions
from $S_{n}$ to $S_{n}$ and $e$ is a constant in $P$. Their interpretation
in $\tilde{M}$ are as follows:
\begin{itemize}
\item $\pi_{1}\left(\left(u,v,i\right)\right)=u$, $\pi_{2}\left(\left(u,v,i\right)\right)=v$. 
\item For each $1\leq n<\omega$, $S_{n}\left(v\right)\Leftrightarrow\left|v\right|\leq n$.
\item For each $1\leq n$, $u\subseteq_{n}v$ if and only if $\left|u\right|\leq n$,
$\left|v\right|\leq n$ and $u\subseteq v$.
\item $u\cap_{n}v=u\cap v$ for all $u,v\in S_{n}$. 
\item $u-_{n}v=u\backslash v$ for $v,u\in S_{n}$. 
\item $u\cup_{n}v=u\cup v$ for $u,v\in S_{n}$. 
\item $e=\emptyset$.
\end{itemize}
Note that they are indeed definable:
\begin{enumerate}
\item $\pi_{1}\left(x\right)$ is the unique $y$ such that $Q_{1}\left(y,x\right)$,
and similarly $\pi_{2}$ is definable.
\item Let $E\left(x,y\right)$ by an auxiliary equivalence relation defined
by $\pi_{1}\left(x\right)=\pi_{1}\left(y\right)\land\pi_{2}\left(x\right)=\pi_{2}\left(y\right)$.
\item $e$ is the unique element $x\in P$ such that there exists exactly
one element $z\in R$ such that $\pi_{1}\left(z\right)=x=\pi_{2}\left(z\right)$. 
\item $x\subseteq_{n}y$ is defined by ``$P\left(x\right),P\left(y\right)$
and the number of elements in the $E$ class of some (equivalently
any) element $z$ such that $\pi_{1}\left(z\right)=x$, $\pi_{2}\left(z\right)=y$
is at most $n+1$''.
\item $S_{n}\left(x\right)$ is defined by ``$P\left(x\right)$ and $e\subseteq_{n}x$''
(In particular, $e\in S_{n}$ for all $n$). 
\item $\cap_{n}$ and $-_{n}$ are then naturally definable using $\subseteq_{n}$.
For instance $x-_{n}y=z$ if and only if $x,y,z$ are in $S_{n}$,
$z\subseteq_{n}x$ and for each $e\neq w\subseteq_{n}y$, $w\nsubseteq_{n}z$. 
\item $x\cup_{n}y=z$ if and only if $x,y\in S_{n}$, $z\in S_{2n}$, $x,y\subseteq_{2n}z$
and $z-_{2n}x\subseteq_{2n}y$. 
\end{enumerate}
Furthermore, $\mathordi{\subseteq_{k}\upharpoonright S_{n}}=\mathordi{\subseteq_{n}}$
for $n\leq k$. Hence every model $M$ of $T$ gives rise naturally
to an induced PBA: $B^{M}:=\left(\bigcup_{n}S_{n}^{M},\cup^{M},\cap^{M},-^{M},\subseteq^{M},e^{M}\right)$
where $\cup^{M}=\bigcup\left\{ \cup_{n}^{M}\left|\, n<\omega\right.\right\} $,
and similarly for $\subseteq^{M},-^{M}$ and $\cap^{M}$ (see Definition
\ref{def:PBA} above).
\begin{claim*}
In the extended language, $T$ eliminates quantifiers.\end{claim*}
\begin{proof}
Suppose $M,N\models T$ are saturated models, $\left|M\right|=\left|N\right|$
and $A\subseteq M,N$ is a common substructure (where $\left|A\right|<\left|M\right|$).
It is enough to show that we have an isomorphism from $M$ to $N$
fixing $A$.

By Proposition \ref{pro:ISomPBA}, we have an isomorphism $f$ from
$B^{M}$ to $B^{N}$ preserving $A\cap B^{M}$ (by saturation and
the choice of language, the condition of the proposition are satisfied).

On $P^{M}\backslash\left(B^{M}\cup A\right)$ there is no structure
and it has the same size as $P^{N}\backslash\left(B^{N}\cup A\right)$
(namely $\left|N\right|$), so we can extend the isomorphism $f$
to $P^{M}$. 

We are left with $R^{M}$: let $a\in R^{M}$, and $a_{i}=\pi_{i}\left(a\right)$
for $i=1,2$. We already defined $f\left(a_{1}\right),f\left(a_{2}\right)$.
Suppose $a_{1}\subseteq_{n}a_{2}$ for minimal $n$. Then there are
exactly $n$ elements $z\in R^{M}$ with $\pi_{1}\left(z\right)=a_{1},\pi_{2}\left(z\right)=a_{2}$.
This is true also in $R^{N}$, and the number of such $z$'s not in
$A$ is the same for both $M$, $N$. Hence we can take this $E$-equivalence
class from $M$ to the appropriate class in $N$.

If not, i.e., $a_{1}\nsubseteq_{n}a_{2}$ for all $n$, then there
are infinitely many elements $z$ in $N$ and in $M$ with $\pi_{1}\left(z\right)=a_{1}$,
$\pi_{2}\left(z\right)=a_{2}$, and again we take this $E$-class
in $M$ outside of $A$ to the appropriate $E$-class in $N$.
\end{proof}
Now we can conclude the proof by a counting types argument. Let $M$
be a countable model of $T$. Let $p\left(x\right)$ be a non-algebraic
type over $M$. There are some cases:
\begin{casenv}
\item $S_{n}\left(x\right)\in p$ for some $n$. Then the type is determined
by the maximal element $c$ in $M$ such that $c\subseteq_{n}x$ (this
is easy, but also follows from the proof of Proposition \ref{pro:ISomPBA}).
\item $S_{n}\left(x\right)\notin p$ for all $n$ but $P\left(x\right)\in p$.
Then $x$ is already determined --- there is nothing more we can say
on $x$. 
\item $R\left(x\right)\in p$. Then the type of $x$ is determined by the
type of $\left(\pi_{1}\left(x\right),\pi_{2}\left(x\right)\right)$
over $M$. 
\end{casenv}
So the number of types over $M$ is countable.\end{proof}
\begin{prop}
Every weakly generic pair of models of $T$ has the independence property.\end{prop}
\begin{proof}
Suppose $\left(M_{1},M\right)$ is a weakly generic pair. We think
of it as a structure of the language $L_{Q}$, where $Q$ is interpreted
as $M$. Consider the formula 
\[
\varphi\left(x,y\right)=P\left(x\right)\land P\left(y\right)\land\exists z\notin Q\left(Q_{1}\left(x,z\right)\land Q_{2}\left(y,z\right)\right).
\]
 This formula has IP: Let $\left\{ a_{i}\left|\, i<\omega\right.\right\} \subseteq M$
be elements from $P^{M}$ such that $a\in S_{1}^{M}$ (as in the language
of the proof of Claim \ref{cla:omegaStable}), i.e., they are atoms
in the induced PBA, and $a_{i}\neq a_{j}$ for $i\neq j$. For any
finite $s\subseteq\omega$ of size $n$, there is an element $b_{s}\in P^{M}$
be such that $a_{i}\subseteq_{n}^{M}b_{s}$ for all $i\in s$. Then
for all $i\in\omega$, $\varphi\left(a_{i},b_{s}\right)$ if and only
if $i\notin s$: 

If $\varphi\left(a_{i},b_{s}\right)$ there are infinitely many $z$'s
in $M$ such that $Q_{1}\left(a_{i},z\right)\wedge Q_{2}\left(b_{s},z\right)$
(otherwise they would all be in $M$). This means that $a_{i}\nsubseteq_{n}^{M}b_{s}$
so $i\notin s$. 

For the other direction, the same exact argument works, but this time
use the fact that the pair is weakly generic.
\end{proof}

\section{\label{sec:Directionality}Directionality}

\subsection{Introduction}
\begin{defn}
A global type $p\left(x\right)\in S\left(\C\right)$ is said to be
\emph{finitely satisfiable} in a set $A$, or a \emph{coheir} over
$A$ if for every formula $\varphi\left(x,y\right)$, if $\varphi\left(x,b\right)\in p$,
then for some $a\in A$, $\varphi\left(a,b\right)$ holds. 
\end{defn}
It is well known (see \cite{Ad}) that a theory $T$ is dependent
if and only if given a type $p\left(x\right)\in S\left(M\right)$
over a model $M$, the number of complete global types $q\in S\left(\C\right)$
that extend $p$ and are finitely satisfiable in $M$ is at most $2^{\left|M\right|}$
(while the maximal number is $2^{2^{\left|M\right|}}$). 

We analyze the behavior of the number of global coheir extensions
in a dependent theory and classify theories by what we call directionality:

Say that $T$ has \emph{small} directionality if and only if the number
of $\Delta$-coheirs (for a finite set of formulas $\Delta$) that
extend a type $p\in S\left(M\right)$ is finite. $T$ has \emph{medium}
directionality if this number is $\left|M\right|$, and it has \emph{large}
directionality if it is neither small or medium. In that case we will
show that it is at least $\ded\left|M\right|$. 

We give an equivalent definition in terms of the number of global
coheir extensions (see Theorem \ref{thm:trichotomy}). 

As far as we know, the first person to give an example of a dependent
theory with large directionality was Delon in \cite{Delon}. 

We give simple combinatorial examples for each of the possible directionalities,
and furthermore we show that RCF and some theories of valued fields
are large.

We do not always assume that $T$ is dependent in this section.

\subsection{Equivalent definitions of directionality}
\begin{defn}
For a type $p\in S\left(A\right)$, let: 
\[
\uf\left(p\right)=\left\{ q\in S\left(\C\right)\left|\, q\mbox{ is a coheir extension of }p\mbox{ over }A\right.\right\} .
\]
For a partial type $p\left(x\right)$ over a set $A$, and a set of
formulas $\Delta$, 
\[
\uf_{\Delta}\left(p\right)=\left\{ q\left(x\right)\in S_{\Delta}\left(\C\right)\left|\, q\cup p\mbox{ is f.s. in }A\right.\right\} .
\]
Note: this definition only makes sense if $p$ is finitely satisfiable
in $A$. The notation $\uf$ refers to ultrafilter.
\end{defn}
And here is the main definition of this section:
\begin{defn}
\label{def:directionality}Let $T$ be any theory, then:
\begin{enumerate}
\item $T$ is said to have \emph{small} directionality (or just, $T$ is
small) if and only if for all finite $\Delta$, $M\models T$ and
$p\in S\left(M\right)$, $\uf_{\Delta}\left(p\right)$ is finite. 
\item $T$ is said to have \emph{medium} directionality (or just, $T$ is
medium) if and only if for every $\lambda\geq\left|T\right|$, 
\[
\lambda=\sup\left\{ \left|\uf_{\Delta}\left(p\right)\right|\left|\, p\in S\left(M\right),\Delta\mbox{ finite},\left|M\right|=\lambda\right.\right\} .
\]
 
\item $T$ is said to have \emph{large} directionality (or just, $T$ is
large) if $T$ is neither small nor medium. 
\end{enumerate}
\end{defn}
\begin{obs}
\label{obs:IPLarge} If $T$ has the independence property, then it
is large. In fact, if $\varphi\left(x,y\right)$ has the independence
property, then there is a type $p\left(x\right)$ over a model $M$,
that has $2^{2^{\left|M\right|}}$ many $\left\{ \varphi\left(x,y\right)\right\} $-extensions
that are finitely satisfiable in $M$. \end{obs}
\begin{proof}
We may assume that $T$ has Skolem functions. Let $\lambda\geq\left|T\right|$,
and let $\bar{a}=\left\langle a_{i}\left|\, i<\lambda\right.\right\rangle $,
$\left\langle b_{s}\left|\, s\subseteq\lambda\right.\right\rangle $
be such that $\left\langle a_{i}\left|\, i<\lambda\right.\right\rangle $
is indiscernible and $\varphi\left(a_{i},b_{s}\right)$ holds iff
$i\in s$. Let $M$ be the Skolem hull of $\bar{a}$. Let $p\left(x\right)\in S\left(M\right)$
be the limit of $\bar{a}$ in $M$ (so $\psi\left(x,c\right)\in p$
iff $\psi\left(x,c\right)$ holds for an end segment of $\bar{a}$).
Let $P\subseteq\SS\left(\lambda\right)$ be an independent family
of size $2^{\lambda}$ (i.e., such that every finite Boolean combination
has size $\lambda$). Then for each $D\subseteq P$, $p\left(x\right)\cup\left\{ \varphi\left(x,b_{s}\right)^{s\in D}\left|\, s\in P\right.\right\} $
is finitely satisfiable in $M$. 
\end{proof}

\subsubsection{Small directionality}

The following construction will be useful (here and in Section \ref{sec:Splintering}):
\begin{const}
\label{const:not small directionality}Let $T$ be any complete theory
and $M\models T$. Suppose that there is some $p\in S\left(M\right)$
and finite $\Delta$ such that $\uf_{\Delta}\left(p\right)$ is infinite,
and contains $\left\{ q_{i}\left|\, i<\omega\right.\right\} $. 

For all $i<j<\omega$, there is a formula $\varphi_{i,j}\in\Delta$
and $b_{i,j}\in\C$ such that $\varphi_{i,j}\left(x,b_{i,j}\right)\in q_{i}$,
$\neg\varphi_{i,j}\left(x,b_{i,j}\right)\in q_{j}$ (or the other
way around). By Ramsey's Theorem we may assume $\varphi_{i,j}$ is
constant --- $\varphi\left(x,y\right)$. Let $N$ be a model containing
$M$ and $\left\{ b_{i,j}\left|\, i<j<\omega\right.\right\} $.

Suppose $\left\langle c_{i}\left|\, i<\omega\right.\right\rangle $
are in $\C$ and $c_{i}\models q_{i}|_{N}$. Let $N'$ be a model
containing $\left\{ c_{i}\left|\, i<\omega\right.\right\} \cup N$.
Let $M^{*}=\left(N',N,M,Q,\bar{f}\right)$ where $Q=\left\{ c_{i}\left|\, i<\omega\right.\right\} $
and $\bar{f}:Q^{2}\to N$ is a tuple of functions of length $\lg\left(y\right)$
defined by $\bar{f}\left(c_{i},c_{j}\right)=\bar{f}\left(c_{j},c_{i}\right)=b_{i,j}$
for $i<j$.

So if $N\models Th\left(M^{*}\right)$ then $N=\left(N_{0}',N_{0},M_{0},Q_{0},\bar{f}_{0}\right)$
and
\begin{itemize}
\item $M_{0}\prec N_{0}\prec N_{0}'\models T$,
\item $N_{0}\cup Q_{0}\subseteq N_{0}'$,
\item $\bar{f}_{0}$ are functions from $Q_{0}^{2}$ to $N_{0}$,
\item For all $c,d\in Q_{0}$, $c\equiv_{M_{0}}d$,
\item $\tp_{\Delta}\left(c/N_{0}\right)\cup\tp\left(c/M_{0}\right)$ is
finitely satisfiable in $M_{0}$ for all $c\in Q_{0}$, and
\item $\varphi\left(c,\bar{f}_{0}\left(c,d\right)\right)\triangle\varphi\left(d,\bar{f}_{0}\left(c,d\right)\right)$
(where $\triangle$ denotes symmetric difference) holds for all $c\neq d\in Q_{0}$.
\end{itemize}
\end{const}
\begin{claim}
\label{cla:bounded}Let $T$ be any theory. Then $T$ is small if
and only if for every $M\models T$ and every type $p\left(x\right)\in S\left(M\right)$,
$\left|\uf\left(p\right)\right|\leq2^{\left|T\right|}$ (here $p$
can also be an infinitary type, but then the bound is $2^{\left|T\right|+\left|\lg\left(x\right)\right|}$). 

In addition, if $T$ is not small, then for every $\lambda\geq\left|T\right|$,
there is a model $M\models T$ of cardinality $\lambda$, a type $p\in S\left(M\right)$,
and a finite set of formulas $\Delta$ such that $\left|\uf_{\Delta}\left(p\right)\right|\geq\lambda$.\end{claim}
\begin{proof}
Assume that $T$ is small. The injective function $\uf\left(p\right)\to\prod_{\varphi\in L}\uf_{\left\{ \varphi\right\} }\left(p\right)$
shows that $\left|\uf\left(p\right)\right|\leq2^{\left|T\right|}$.

Conversely (and the ``In addition'' part): Assume that there is
some $p$ and $\Delta$ such that $\uf_{\Delta}\left(p\right)$ is
infinite. Use Construction \ref{const:not small directionality}: 

For every $\lambda\geq\left|T\right|$ we may find $N\models Th\left(M^{*}\right)$
of size $\lambda$ such that $\left|M_{0}\right|=\left|Q_{0}\right|=\lambda$,
and we have a model $M_{0}$ of $T$ with a type $p$ over it, which
has at least $\lambda$ many $\Delta$-coheirs.
\end{proof}
We conclude this section with a claim on theories with non-small directionality.
\begin{claim}
Suppose $T$ has medium or large directionality. Then there exists
some $M\models T$, $p\in S\left(M\right)$, $\psi\left(x,y\right)$
and $\left\{ c_{i}\left|\, i<\omega\right.\right\} \subseteq\C$ such
that for each $i<\omega$ the set $p\left(x\right)\cup\left\{ \psi\left(x,c_{j}\right)^{j=i}\left|\, j<\omega\right.\right\} $
is finitely satisfiable in $M$. \end{claim}
\begin{proof}
We consider the structure $M^{*}$ introduced in Construction \ref{const:not small directionality}
and the formula $\varphi$ chosen there. Find an extension of $M^{*}$
with an indiscernible sequence $\left\langle d_{i}\left|\, i\in\mathbb{Z}\right.\right\rangle $
inside $Q$. Assume without loss that $\varphi\left(d_{0},\bar{f}\left(d_{0},d_{1}\right)\right)\land\neg\varphi\left(d_{1},\bar{f}\left(d_{0},d_{1}\right)\right)$
holds. This means that $\varphi\left(d_{0},\bar{f}\left(d_{0},d_{1}\right)\right)\land\neg\varphi\left(d_{0},\bar{f}\left(d_{-1},d_{0}\right)\right)$.

We claim that $\varphi\left(d_{i},\bar{f}\left(d_{j},d_{j+1}\right)\right)\land\neg\varphi\left(d_{i},\bar{f}\left(d_{j-1},d_{j}\right)\right)$
holds if and only if $i=j$: 

Suppose this holds but $i\neq j$. If $i>j$ then, since $\neg\varphi\left(d_{1},f\left(d_{0},d_{1}\right)\right)$,
it must be that $i>j+1$, but then we have a contradiction to indiscernibility.
Similarly, it cannot be that $i<j$. Thus the claim is proved with
$\psi\left(x;y,z\right)=\varphi\left(x,y\right)\land\neg\varphi\left(x,z\right)$,
and $c_{i}=\left\langle f\left(d_{i},d_{i+1}\right),f\left(d_{i-1},d_{i}\right)\right\rangle $. 
\end{proof}

\subsubsection{Some helpful facts about dependent theories}

Assume $T$ is dependent.

Recall,
\begin{defn}
A global type $p\left(x\right)$ is \emph{invariant} \emph{over} a
set $A$ if it does not split over it, namely if whenever $b$ and
$c$ have the same type over $A$, $\varphi\left(x,b\right)\in p$
if and only if $\varphi\left(x,c\right)\in p$ for every formula $\varphi\left(x,y\right)$.
\end{defn}

\begin{defn}
\label{def:tensor}Suppose $p\left(x\right)$ and $q\left(y\right)$
are global $A$-invariant types. Then $\left(q\otimes p\right)\left(x,y\right)$
is a global invariant type defined as follows: for any $B\supseteq A$,
let $a_{B}\models p|_{B}$ and $b_{B}\models q|_{Ba_{B}}$, then $p\otimes q=\bigcup_{B\supseteq A}\tp\left(a_{B},b_{B}/B\right)$.
One can easily check that it is well defined and $A$-invariant. Let
$p^{\left(n\right)}=p\otimes p\cdots\otimes p$ where the product
is done $n$ times. So $p^{\left(n\right)}$ is a type in $\left(x_{0},\ldots,x_{n-1}\right)$,
and $p^{\left(\omega\right)}=\bigcup_{n<\omega}p^{\left(n\right)}$
is a type in $\left(x_{0},\ldots,x_{n},\ldots\right)$. For $n\leq\omega$,
$p^{\left(n\right)}$ is a type of an $A$-indiscernible sequence
of length $n$. \end{defn}
\begin{fact}
\emph{\label{fac:pOmega}}\cite[Lemma 2.5]{HP} If $T$ is NIP then
for a set $A$ the map $p\left(x_{0}\right)\mapsto p^{\left(\omega\right)}\left(x_{0},\ldots\right)|_{A}$
from global $A$-invariant types to $\omega$-types over $A$ is injective.
\end{fact}
In the rest of the section, $\Delta$ will always denote a finite
set of formulas, closed under negation.
\begin{claim}
\label{cla:local types define indiscernibles} For every set $A\subseteq C$,
any type $q\left(x\right)\in S_{\Delta}\left(\C\right)$ which is
finitely satisfiable in $A$ and any choice of a coheir $q'\in S\left(\C\right)$
over $A$ which completes $q$:
\begin{itemize}
\item $\left(a_{0},\ldots,a_{n-1}\right)\models\left(q'^{\left(n\right)}|_{C}\right)\upharpoonright\Delta$
if and only if $a_{0}\models q|_{C}$, $a_{1}\models q|_{Ca_{0}}$,
etc. 
\end{itemize}
This enables us to define $q^{\left(n\right)}\left(x_{0},\ldots x_{n-1}\right)\in S_{\Delta}\left(\C\right)$
as $q'^{\left(n\right)}\upharpoonright\Delta$. 

It follows that $q^{\left(n\right)}$ is a type of a $\Delta$-indiscernible
sequence of length $n$. \end{claim}
\begin{proof}
The proof is by induction on $n$: 

Right to left: suppose $a_{i}\models q|_{Ca_{0}\cdots a_{i-1}}$ for
$i\leq n$, $\varphi\left(x_{0},\ldots,x_{n},y\right)\in\Delta$ and
$\varphi\left(x_{0},\ldots,x_{n},c\right)\in q'^{\left(n+1\right)}$
for $c\in C$ but $\neg\varphi\left(a_{0},\ldots,a_{n},c\right)$
holds. Then by the choice of $a_{n}$, $\neg\varphi\left(a_{0},\ldots,a_{n-1},x,c\right)\in q$.
Suppose $\left(b_{0},\ldots,b_{n-1}\right)\models q'^{\left(n\right)}|_{C}$,
then $\varphi\left(b_{0},\ldots,b_{n-1},x,c\right)\in q$ so there
is some $c'\in A$ such that $\varphi\left(b_{0},\ldots,b_{n-1},c',c\right)\land\neg\varphi\left(a_{0},\ldots,a_{n-1},c',c\right)$
holds. But this is a contradiction to the induction hypothesis.

Left to right is similar. 
\end{proof}
The following is a local version of Fact \ref{fac:pOmega}, which
will be useful later:
\begin{prop}
\label{prop:localpomega}($T$ dependent) Suppose $\Delta$ is a finite
set of formulas, $x$ a finite tuple of variables. Then there exists
$n<\omega$ and finite set of formulas $\Delta_{0}$ such that for
every set $A$, if $q_{1}\left(x\right),q_{2}\left(x\right)\in S\left(\C\right)$
are coheirs over $A$ and $\left(q_{1}^{\left(n\right)}\upharpoonright\Delta_{0}\right)|_{A}=\left(q_{2}^{\left(n\right)}\upharpoonright\Delta_{0}\right)|_{A}$
then $q_{1}\upharpoonright\Delta=q_{2}\upharpoonright\Delta$.\end{prop}
\begin{proof}
By compactness and NIP,
\begin{itemize}
\item there exists some finite set of formulas $\Delta_{0}$ and some $n$
such that for all $\varphi\left(x,y\right)\in\Delta$ and all $\Delta_{0}$-indiscernible
sequences $\left\langle a_{0},\ldots,a_{n_{\Delta}-1}\right\rangle $,
there is \uline{no} $c$ such that $\varphi\left(a_{i},c\right)$
holds if and only if $i$ is even. We may assume that $\Delta\subseteq\Delta_{0}$. 
\end{itemize}
By Claim \ref{cla:local types define indiscernibles}, we can conclude:
\begin{itemize}
\item Suppose that $\left(q_{1}^{\left(n\right)}\upharpoonright\Delta_{0}\right)|_{A}=\left(q_{2}^{\left(n\right)}\upharpoonright\Delta_{0}\right)|_{A}$,
but $q_{1}\upharpoonright\Delta\neq q_{2}\upharpoonright\Delta$.
Then there is some formula $\varphi\left(x,y\right)\in\Delta$ and
some $c\in\C$ such that $\varphi\left(x,c\right)\in q_{1}$ and $\neg\varphi\left(x,c\right)\in q_{2}$.

Since $\left(q_{1}\upharpoonright\Delta_{0}\right)^{\left(n\right)}|_{A}=\left(q_{2}\upharpoonright\Delta_{0}\right)^{\left(n\right)}|_{A}$,
$\left(q_{1}\upharpoonright\Delta_{0}\right)^{\left(m\right)}|_{A}=\left(q_{2}\upharpoonright\Delta_{0}\right)^{\left(m\right)}|_{A}$
for every $m\leq n$, and it follows by induction on $m$ that the
sequence defined by $a_{0}\models\left(q_{1}\upharpoonright\Delta_{0}\right)|_{Ac}$,
$a_{1}\models\left(q_{2}\upharpoonright\Delta_{0}\right)|_{Aca_{0}}$,
$a_{2}\models\left(q_{1}\upharpoonright\Delta_{0}\right)|_{Aca_{0}a_{1}}$,
$\ldots$, $a_{m-1}\models\left(q_{i}\upharpoonright\Delta_{0}\right)|_{Aca_{0}\cdots a_{m-2}}$
($i\in\left\{ 1,2\right\} $) realizes this type. But this entails
a contradiction, because $\left\langle a_{0},\ldots,a_{n-1}\right\rangle $
is a $\Delta_{0}$ indiscernible sequence (even over $A$), while
$\varphi\left(a_{i},c\right)$ holds if and only if $i$ is even.

\end{itemize}
\end{proof}
\begin{problem}
Does Proposition \ref{prop:localpomega} hold for invariant types
(not just for coheirs)?
\end{problem}

\subsubsection{Large directionality and definability }

Let us recall the definition of $\ded\lambda$. 
\begin{defn}
\label{def:ded}Let $\ded\lambda$ be the supremum of the set: 
\[
\left\{ \left|I\right|\left|\, I\mbox{ is a linear order with a dense subset of size }\leq\lambda\right.\right\} .
\]
\end{defn}
\begin{fact}
\label{fac:Ded}It is well known that $\lambda<\ded\lambda\leq\left(\ded\lambda\right)^{\aleph_{0}}\leq2^{\lambda}$.
If $\lambda^{<\lambda}=\lambda$ then $\ded\lambda=2^{\lambda}$ so
$\ded\lambda=\left(\ded\lambda\right)^{\aleph_{0}}=2^{\lambda}$. 
\end{fact}
For more, see Section \ref{sec:Appendix:-dense-types} and \cite[Section 6]{Sh1007}.
\begin{defn}
\label{def:ExtDefType}Suppose $M$ is a model and $p\in S\left(M\right)$.
Let $M_{p}$ be $M$ enriched with externally definable sets defined
over a realization of $p$. Namely, we enrich the language to a language
$L_{p}$ by adding new relation symbols $\left\{ d_{p}x\varphi\left(x,y\right)\left|\,\varphi\left(x,y\right)\mbox{ is a formula}\right.\right\} $
(so $d_{p}$ is thought of as a quantifier over $x$), and let $M_{p}$
be a structure for $L_{p}$ with universe $M$ where we interpret
$d_{p}x\varphi\left(x,y\right)$ as $\left\{ b\in M\left|\,\varphi\left(x,b\right)\in p\right.\right\} $. \end{defn}
\begin{rem}
\label{rem:TpModel}Every model $N\models Th\left(M_{p}\right)$ gives
rise to a complete $L$ type over $N$, namely $p^{N}=\left\{ \varphi\left(b,x\right)\left|\, b\in N,N\models d_{p}x\varphi\left(x,b\right)\right.\right\} $.\end{rem}
\begin{claim}
\label{cla:Large}Let $T$ be any theory, $M\models T$. Suppose $p\in S\left(M\right)$,
$q\in\uf\left(p\right)$, and $\bar{a}=\left\langle a_{0},a_{1},\ldots\right\rangle \models q^{\left(\omega\right)}|_{M}$.
If $\tp\left(\bar{a}/M\right)$ is not definable with parameters in
$M_{p}$, then $T$ is large. 

Moreover, in this case 
\begin{itemize}
\item [$\otimes$]There exists a finite $\Delta$ such that for every $\lambda\geq\left|T\right|$,
\[
\ded\lambda\leq\sup\left\{ \left|\uf_{\Delta}\left(p\right)\right|\left|\, p\in S\left(N\right),\left|N\right|=\lambda\right.\right\} .
\]

\end{itemize}
\end{claim}
\begin{proof}
We may assume that $\left|M\right|=\left|L\right|$: let $r=q^{\left(\omega\right)}|_{M}$,
and $N\prec M_{r}$, $\left|N\right|=\left|L\right|$. Then $N$ gives
rise to a complete type $r'\left(x_{0},x_{1},\ldots\right)\in S\left(N\right)$.
Let $p'=r'\upharpoonright x_{0}$. It is easy to see that $r'=q'^{\left(\omega\right)}|_{N}$
for some $q'\in\uf\left(p'\right)$. Also, $r'$ is not definable
with parameters in $N_{p'}$. 

Let us recall a theorem from \cite{Sh009} (we formulate it a bit
differently):

Suppose $L$ is a language of cardinality at most $\lambda$, $P$
a new predicate (or relation symbol), and $S$ a complete theory in
$L\left(P\right)$.
\begin{defn*}
$\DfOne\lambda$ is the the supremum of the set of cardinalities:
\[
\left|\left\{ B'\subseteq M\left|\,\left(M,B'\right)\cong\left(M,B\right)\right.\right\} \right|
\]
where $\left(M,B\right)$ is an $L\left(P\right)$ model of $S$ of
cardinality $\lambda$. \end{defn*}
\begin{thm*}
\cite[Theorem 12.4.1]{Sh009,Hod} The following are equivalent%
\footnote{The original theorem referred to $\ded^{*}\lambda$, which counts
the number of branches of the same height in a tree with $\lambda$
many nodes, but it equals $\ded\lambda$, see \cite[Section 6]{Sh1007}
and Fact \ref{fac:calculating ded}.%
}:
\begin{enumerate}
\item $P$ is not definable with parameters in $S$, i.e., there is no $L$-formula
$\theta\left(x,y\right)$ such that $S\models\exists y\forall x\left(P\left(x\right)\leftrightarrow\theta\left(x,y\right)\right)$.
\item For every $\lambda\geq\left|L\right|$, $\DfOne\left(\lambda\right)\geq\ded\lambda$.
\end{enumerate}
\end{thm*}
Let $n<\omega$ be the integer first such that $\tp\left(\bar{a}\upharpoonright n/M\right)$
is not definable with parameters in $M_{p}$. So $1<n$ and $r=\tp\left(a_{0},\ldots,a_{n-2}/M\right)$
is definable but $\tp\left(a_{0},\ldots,a_{n-1}/M\right)$ is not. 

For a formula $\alpha\left(x_{0},\ldots,x_{n-2},y\right)$ let $\left(d_{r}\alpha\right)\left(y\right)$
be a formula in $L\left(M_{p}\right)$ defining $\alpha\left(\bar{a}\upharpoonright n-1,M\right)$.
If $M_{p}\prec N\models T_{p}$ then, as in Remark \ref{rem:TpModel},
there is a complete $L$ type $r^{N}\left(x_{0},\ldots,x_{n-2}\right)$
over $N$ defined by $\alpha\left(x_{0},\ldots,x_{n-2},b\right)\in r^{N}$
if and only if $N\models\left(d_{r}\alpha\right)\left(b\right)$.

There is some formula $\varphi\left(x_{0},\ldots,x_{n-1},y\right)$
such that the set $B_{0}:=\varphi\left(a_{0},\ldots,a_{n-1},M\right)$
is not definable with parameters in $M_{p}$. Let $S=Th\left(M_{p},B_{0}\right)$
in the language $L\left(M_{p}\right)\left(P\right)$ (naming elements
from $M$, so that $N\models S$ implies $M_{p}\prec N$).

By the theorem cited above, for every $\lambda\geq\left|L\right|$
and $\kappa<\ded\lambda$, there exists a model $\left(N_{\lambda,\kappa},B_{\lambda,\kappa}\right)=\left(N,B\right)\models S$
of cardinality $\lambda$ such that, letting $\mathcal{B}^{N}=\left\{ B'\left|\,\left(N,B'\right)\cong\left(N,B\right)\right.\right\} $,
$\left|\mathcal{B}^{N}\right|>\kappa$.

Let $\bar{a}^{N}=\left(a_{0}^{N},\ldots,a_{n-2}^{N}\right)\models r^{N}$
and for every $B'\in\mathcal{B}^{N}$, let 
\[
q_{B'}=p^{N}\left(x\right)\cup\left\{ \varphi\left(\bar{a}^{N},x,\bar{b}\right)\left|\,\bar{b}\in B'\right.\right\} \cup\left\{ \neg\varphi\left(\bar{a}^{N},x,\bar{b}\right)\left|\,\bar{b}\notin B'\right.\right\} .
\]
By choice of $S$, $B$ and $\mathcal{B}$, $q_{B'}$ is finitely
satisfiable in $N$ and for $B'\neq B''\in\mathcal{B}^{N}$, $q_{B'}\upharpoonright\varphi\neq q_{B''}\upharpoonright\varphi$,
so now $N\upharpoonright L$ is a model of $T$ with a type $p^{N}$
such that $\left|\uf_{\left\{ \varphi\right\} }\left(p\right)\right|>\kappa$. 
\end{proof}
If $T$ is small we can say more: 
\begin{claim}
Assume $T$ is small, $M\models T$. Suppose $p\in S\left(M\right)$,
$q\in\uf\left(p\right)$, and $\bar{a}=\left\langle a_{0},a_{1},\ldots\right\rangle \models q^{\left(\omega\right)}|_{M}$,
then $\tp\left(\bar{a}/M\right)$ is definable over $\acl^{\eq}\left(\emptyset\right)$
in $M_{p}$.\end{claim}
\begin{proof}
By Claim \ref{cla:Large} it is definable with parameters in $M_{p}$.
Let $n<\omega$ be minimal such that $q^{\left(n\right)}|_{M}$ is
not definable over $\acl^{\eq}\left(\emptyset\right)$ in $M_{p}$.
Suppose that for some formula $\varphi\left(x_{0},\ldots,x_{n-1},y\right)$,
$\tp_{\varphi}\left(\bar{a}/M\right)$ is not definable over $\acl^{\eq}\left(\emptyset\right)$
in $M_{p}$. This means that while the set $\left\{ b\in M\left|\,\C\models\varphi\left(\bar{a},b\right)\right.\right\} $
is definable by $\psi\left(y,c\right)$ for some $\psi$ in $L\left(M_{p}\right)$,
$c\notin\acl^{\eq}\left(\emptyset\right)$. We may assume that $c\in M_{p}^{\eq}$
is the code of this set (for every automorphism $\sigma$ of the monster
model $\C$ of $M_{p}^{\eq}$, $\sigma$ fixes $\psi\left(\C,c\right)$
if and only if $\sigma\left(c\right)=c$). So in some elementary extension
$M_{p}\prec N$, there are infinitely many conjugates of $c$ over
$\acl^{\eq}\left(\emptyset\right)$, $\left\{ c_{i}\left|\, i<\omega\right.\right\} $,
such that $\psi\left(N,c_{i}\right)\neq\psi\left(N,c_{j}\right)$
for $i\neq j$. This implies that $\uf_{\left\{ \varphi\right\} }\left(p^{N}\right)\geq\aleph_{0}$,
just as in the proof of Claim \ref{cla:Large}.\end{proof}
\begin{cor}
\label{cor:Medium}($T$ dependent) $T$ is large if and only if for
every $\lambda\geq\left|T\right|$, 
\[
\sup\left\{ \left|\uf_{\Delta}\left(p\right)\right|\left|\, p\in S\left(M\right),\Delta\mbox{ finite},\left|M\right|=\lambda\right.\right\} =\ded\lambda.
\]
\end{cor}
\begin{proof}
Suppose $T$ is large, i.e., for some $\Delta$, $M$ of size $\lambda\geq\left|T\right|$
and $p\in S\left(M\right)$, $\left|\uf_{\Delta}\left(p\right)\right|>\lambda$.
By Proposition \ref{prop:localpomega}, find some $\Delta_{0}$ and
$n$ such that 
\[
\left|\uf_{\Delta}\left(p\right)\right|\leq\left|\left\{ \left(q^{\left(n\right)}\upharpoonright\Delta_{0}\right)|_{M}\left|\, q\in\uf\left(p\right)\right.\right\} \right|.
\]
Hence there is some $q\in\uf_{\Delta}\left(p\right)$ such that $q^{\left(n\right)}\upharpoonright\Delta_{0}$
is not definable with parameters over $M_{p}$, and we are done by
Claim \ref{cla:Large} (also note that $\left|S_{\Delta}\left(M\right)\right|\leq\ded\left|M\right|$
for every finite $\Delta$ in dependent theories (see e.g., \cite[Theorem 4.3]{Sh10})). 
\end{proof}

\subsubsection{Concluding remarks }
\begin{thm}
\label{thm:trichotomy}For every theory $T$, 
\begin{enumerate}
\item $T$ is small iff for all $M\models T$, $p\left(x\right)\in S\left(M\right)$,
$\left|\uf\left(p\right)\right|\leq2^{\left|T\right|+\left|\lg\left(x\right)\right|}$
.
\item $T$ is medium iff for all $M\models T$, $p\left(x\right)\in S\left(M\right)$,
$\left|\uf\left(p\right)\right|\leq\left|M\right|^{\left|T\right|+\left|\lg\left(x\right)\right|}$
and $T$ is not small.
\end{enumerate}
\end{thm}
\begin{proof}
(1) is Claim \ref{cla:bounded}. 

(2) Left to right is clear. Conversely, since $T$ is not small, by
Claim \ref{cla:bounded}, for all $\lambda\geq\left|T\right|$, 
\[
\sup\left\{ \left|\uf_{\Delta}\left(p\right)\right|\left|\, p\in S\left(M\right),\Delta\mbox{ finite},\left|M\right|=\lambda\right.\right\} \geq\lambda.
\]
 On the other hand, if it is strictly greater than $\lambda$, then
by definition $T$ is large. For $\lambda$ large enough, $\ded\lambda>\lambda=\lambda^{\left|T\right|+\left|\lg\left(x\right)\right|}$
so by Corollary \ref{cor:Medium}, we get a contradiction to the right
hand side. 
\end{proof}
In Section \ref{sub:examples of directionality}, we will show that
these classes are not empty, and thus:
\begin{cor}
For $\lambda\geq\left|T\right|$, the cardinality: 
\[
\sup\left\{ \left|\uf_{\Delta}\left(p\right)\right|\left|\, p\in S\left(M\right),\Delta\mbox{ finite},\left|M\right|=\lambda\right.\right\} 
\]
 has four possibilities: finite / $\aleph_{0}$ ; $\lambda$; $\ded\lambda$;
$2^{2^{\lambda}}$.

This corresponds to small, medium, and large directionality (the last
one happens when the theory has the independence property, see Observation
\ref{obs:IPLarge}). \end{cor}
\begin{problem}
Suppose $T$ is interpretable in $T'$, and $T$ is large. Does this
imply that $T'$ is large or at least not small?
\end{problem}

\subsection{\label{sub:examples of directionality}Examples of different directionalities.}

Here we give examples of the different directionalities.

\subsubsection{Small directionality}
\begin{example}
\label{exa:DLO} $Th\left(\mathbb{Q},<\right)$ has small directionality.
In fact, every $1$-type over a model $M$, has at most 2 global coheirs,
and in general, a type $p\left(x_{0},\ldots,x_{n-1}\right)$ is determined
by the order type of $\left\{ x_{0},\ldots,x_{n-1}\right\} $ and
$p\upharpoonright x_{0}$, $p\upharpoonright x_{1}$, etc. \end{example}
\begin{prop}
\label{pro:DTRBounded}The theory of dense trees is also small. \end{prop}
\begin{proof}
So here $T$ is the model completion of the theory of trees in the
language $\left\{ <,\wedge\right\} $. 
\begin{claim*}
Let $M\models T$ and $p\left(x_{0},x_{1},\ldots,x_{n-1}\right)\in S\left(M\right)$
be any type. Then $\bigcup_{i,j<n}p\upharpoonright\left(x_{i},x_{j}\right)\cup p|_{\emptyset}\vdash p$. \end{claim*}
\begin{proof}
Let $\Sigma=\bigcup_{i,j<n}p\upharpoonright\left(x_{i},x_{j}\right)\cup p|_{\emptyset}$.
Suppose $\left(a_{0},\ldots,a_{n-1}\right)\models\Sigma$. By quantifier
elimination, the formulas in $p$ are Boolean combination of formulas
of the form $\bigwedge_{k<m}x_{j_{k}}\wedge a\leq\bigwedge_{l<r}x_{r_{l}}\wedge b$
where $b,a\in M$ and $j_{0},\ldots,j_{k-1},r_{0},\ldots,r_{l-1}<n$.

If $a,b$ does not appear, $\bar{a}$ satisfy this formula because
we included $p|_{\emptyset}$. Consider $\bigwedge_{k<m}x_{j_{k}}\wedge a$:
by assumption we know what is the ordering of $\left\{ x_{j_{k}}\wedge a\left|\, k<m\right.\right\} $
(this set is linearly ordered --- it is below $a$). Hence, as $\bar{a}\models\Sigma$,
$\bigwedge_{k<m}a_{j_{k}}\wedge a$ must be equal to the minimal element
in this set, namely $a_{j_{k}}\wedge a$ for some $k$, which is determined
by $\Sigma$. Now $a_{j_{k}}\wedge a\leq\bigwedge_{l<r}a_{r_{l}}\wedge b$
holds if and only if for each $r<l$, we have $a_{j_{k}}\wedge a\leq a_{r}$
and $a_{j_{k}}\wedge a\leq b$, both decided in $\Sigma$. 

Note that we can get rid of $p|_{\emptyset}$ but we should replace
2-types by 3-types.\end{proof}
\begin{claim*}
For any $A\subseteq M\models T$, and $p\left(x\right),q\left(y\right)\in S\left(A\right)$,
there are only finitely many complete type $r\left(x,y\right)$ that
contain both $p$ and $q$. In fact there is a uniform bound on their
number.\end{claim*}
\begin{proof}
We may assume that $A$ is a substructure. For any $a\in M$, the
structure generated by $A$ and $a$, denoted by $A\left(a\right)$,
is just $A\cup\left\{ a_{A},a\right\} $ where $a_{A}=\max\left\{ b\wedge a\left|\, b\in A\right.\right\} $.
Note that $a_{A}$ need not exist, but if it does, then it is the
only new element apart from $a$ (because if $b_{1}\wedge a<b_{2}\wedge a$
then $b_{1}\wedge a=b_{1}\wedge b_{2}$).

Now, let $a\models p$ and $b\models q$. Let $B=B\left(p,q\right)=\left\{ d\in A\left|\, d\leq a\,\&\, d\leq b\right.\right\} $.
This set is linearly ordered, and it may have a maximum. If it does,
denote it by $m$. Note that $m$ depends only on $p$ and $q$.

Now it is easy to show that $\tp\left(a,a_{A},b,b_{A},a\wedge b/m\right)\cup\tp\left(a/A\right)\cup\tp\left(b/A\right)$
determines $\tp\left(a,b/A\right)$ by quantifier elimination. This
suffices because the number of types of finite tuples over a finite
set is finite.
\end{proof}
Let $M\models T$, $p\left(\bar{x}\right)\in S\left(M\right)$ and
$I$ be the set of all types $r\left(\bar{x}_{0},\bar{x}_{1},\ldots\right)$
over $M$ such that realizations of $r$ are indiscernible sequences
of tuples satisfying $p$.

Let $V=\left\{ \left(y_{0},y_{1}\right)\left|\, y_{0},y_{1}\mbox{ are any 2 variables from }\bar{x}_{0},\bar{x}_{1}\right.\right\} $.
By the second claim, for any $\left(y_{0},y_{1}\right)\in V$, the
set $\left\{ r\upharpoonright\left(y_{0},y_{1}\right)\left|\, r\in I\right.\right\} $
is finite. By the first claim and indiscernibility, the function taking
$r$ to $\left(r|_{\emptyset},\left\langle r\upharpoonright\left(y_{0},y_{1}\right)\left|\,\left(y_{0},y_{1}\right)\in V\right.\right\rangle \right)$
is injective. Together, it means that $\left|I\right|\leq2^{\aleph_{0}+\left|\lg\left(\bar{x}\right)\right|}$
and we are done by Fact \ref{fac:pOmega}. 
\end{proof}

\subsubsection{Medium directionality}
\begin{example}
Let $L=\left\{ P,Q,H,<\right\} $ where $P$ and $Q$ are unary predicates,
$H$ is a unary function symbol and $<$ is a binary relation symbol.
Let $T^{\forall}$ be the following theory:
\begin{itemize}
\item $P\cap Q=\emptyset$.
\item $H$ is a function from $P$ to $Q$ (so $H\upharpoonright Q=\id$). 
\item $\left(P,<,\wedge\right)$ is a tree. 
\end{itemize}
And let $T$ be its model completion (so $T$ eliminates quantifiers).
Note that there is no structure on $Q$. So as in Section \ref{sec:Indisc},
$T$ is dependent (this theory is interpretable in the theory there). 

Let $T'$ be the restriction of $T$ to the language $L'=L\backslash\left\{ H\right\} $.
The same ``moreover'' part applies here as in Corollary \ref{cor:ModelCom},
so $T'$ is the model completion of $T^{\forall}\upharpoonright L'$
and also eliminates quantifiers. \end{example}
\begin{claim}
$T'$ has small directionality.\end{claim}
\begin{proof}
The only difference between $T'$ and dense trees is the new set $Q$
which has no structure. Easily this does not make any difference. \end{proof}
\begin{prop}
\label{pro:MediumExample}$T$ has medium directionality.\end{prop}
\begin{proof}
Let $M\models T$. Let $B$ be a branch in $P^{M}$ (i.e., a maximal
linearly ordered set). Let $p\left(x\right)\in S\left(M\right)$ be
a complete type containing $\Sigma_{B}:=\left\{ b<x\left|\, b\in B\right.\right\} $.
Note that $\Sigma$ ``almost'' isolates $p$: the only freedom we
have, is to determine what is $H\left(x\right)$. So suppose $H\left(x\right)=m\in p$
for $m\in Q$.

Let $c\models p$ (so $c\notin M$). For each $a\in Q^{M}$, let $p_{a}\left(x\right)=p\cup\left\{ H\left(c\wedge x\right)=a\right\} $.

Then $p_{a}$ is finitely satisfiable in $M$: Suppose $\Gamma\subseteq p_{a}$
is finite. By quantifier elimination, we may assume that $\Gamma\subseteq\Sigma_{B}\cup\left\{ H\left(c\wedge x\right)=a\right\} \cup\left\{ H\left(x\right)=m\right\} $.
Since $B$ is linearly ordered, we may assume that $\Gamma=\left\{ b<x,H\left(c\wedge x\right)=a,H\left(x\right)=m\right\} $
for some $b\in B$. Since $T$ is the model completion of $T^{\forall}$
which has the amalgamation property, there are two elements $d,e\in M$
such that $b<d,e$, $e\in B$, $d\notin B$, $H\left(d\right)=m$
and $H\left(d\wedge e\right)=a$. Since $d\wedge c=d\wedge e$, $d\models\Gamma$.
We have found $\left|Q^{M}\right|$ coheirs of $p$, and since $M$
was arbitrary $T$ is not small. 

This gives a lower bound on the directionality of $T$, and we would
like to find an upper bound as well. We shall use the same idea as
in the proof of Proposition \ref{pro:DTRBounded}. 

Let $M\models T$, $p\left(\bar{x}\right)\in S\left(M\right)$ and
$I$ be the set of all types $r\left(\bar{x}_{0},\bar{x}_{1},\ldots\right)$
over $M$ such that realizations of $r$ are indiscernible sequences
of tuples satisfying $p$. Let $I'=\left\{ r\upharpoonright L'\left|\, r\in I\right.\right\} $.
By the proof of Proposition \ref{pro:DTRBounded}, $\left|I'\right|\leq2^{\aleph_{0}+\left|\lg\left(\bar{x}\right)\right|}$.

Let $V=\left\{ t\left|\, t\mbox{ is a term in }L\mbox{ in the variables }\bar{x}_{0},\bar{x}_{1},\ldots\mbox{ over }\emptyset\right.\right\} $.
Suppose $r\in I$, then, as in the proof of Proposition \ref{pro:DTRBounded},
for every term $t$ such that $P\left(t\right)\in r$, let $t_{r}=\max\left\{ a\wedge t\left|\, a\in M\right.\right\} $
--- a term over $M$ (it need not exist). To determine $r$, it is
enough to determine the equations that occur between the images under
$H$ of the $t_{r}$'s and the $t$'s over $M$. This shows that $\left|I\right|\leq\left|M\right|^{\aleph_{0}+\left|\lg\left(\bar{x}\right)\right|}$. 
\end{proof}

\subsubsection{Large directionality}
\begin{example}
\label{exa:Large}Let $L=\left\{ P,Q,H,<_{P},<_{Q}\right\} $ where
$P$ and $Q$ are unary predicates, $H$ is a unary function symbol
and $<_{P},<_{Q}$ are binary relation symbols. Let $T^{\forall}$
be the following theory:
\begin{itemize}
\item $P\cap Q=\emptyset$.
\item $H$ is a function from $P$ to $Q$.
\item $\left(P,<_{P},\wedge\right)$ is a tree.
\item $\left(Q,<_{Q}\right)$ is a linear order.
\end{itemize}
\end{example}
\begin{prop}
$T$ has large directionality.\end{prop}
\begin{proof}
This is similar to the proof of Proposition \ref{pro:MediumExample}.

Let $M\models T$. Let $B$ be a branch in $P^{M}$. Let $p\left(x\right)\in S\left(M\right)$
be a complete type containing $\Sigma_{B}:=\left\{ b<x\left|\, b\in B\right.\right\} $
saying that $H\left(x\right)=m$ for some $m\in Q^{M}$.

Let $c\models p$ (so $c\notin M$). For each cut $I\subseteq Q^{M}$,
let: 
\[
p_{I}\left(x\right)=p\cup\left\{ e<H\left(c\wedge x\right)<f\left|\, e\in I,f\in Q^{M}\backslash I\right.\right\} .
\]
Then $p_{I}$ is finitely satisfiable in $M$ as in the proof of \ref{pro:MediumExample}.
So for every cut in $Q$ we found a coheir of $p$, and since $M$
was arbitrary $T$ is not small nor medium (because for every linear
order, we can find a model such that $Q$ contains this order). 
\end{proof}

\subsubsection{RCF}

It turns out that even RCF has large directionality, as we shall present
now. 

Apparently, that RCF  was not small was already known and can be deduced
from Marcus Tressl's thesis (see \cite[18.13]{Tressle}), but here
we give a direct proof that RCF is in fact large and even more. 
\begin{defn}
\label{def:dense types}Let $M\models RCF$. A type $p\in S\left(M\right)$
is called \emph{dense} if it is not definable and the differences
$b-a$ with $a,b\in M$ and $a<x<b\in p$, are arbitrarily (w.r.t.
$M$) close to $0$. 
\end{defn}
For example, if $R$ is the real closure of $\mathbb{Q}$, then $\tp\left(\pi/R\right)$
is dense. 
\begin{fact}
\label{fac:DenseExist}Any real closed field can be embedded into
a real closed field of the same cardinality with some dense type.\end{fact}
\begin{proof}
{[}due to Marcus Tressl{]} Let $R$ be a real closed field. Let $S$
be the (real closed) field $R\left(\left(t^{\mathbb{Q}}\right)\right)$
of generalized power series over $R$. Let $K$ be the definable closure
of $R\left(t\right)$ in $S$ and let $p$ be the 1-type of the formal
Taylor series of $e^{t}$ over $K$: $\tp\left(1+t^{1}/1!+t^{2}/2!+\cdots/K\right)$
. Then $p$ is a dense 1-type over $K$. \end{proof}
\begin{claim}
\label{cla:WeakOrth}Suppose $p$ is dense and $q$ is a definable
type over $M\models RCF$ and both are complete. Then $q$ and $p$
are weakly orthogonal, meaning that $p\left(x\right)\cup q\left(y\right)$
implies a complete type over $M$.\end{claim}
\begin{proof}
{[}Remark: this is an easy result that is well known, but for completeness
we give a proof.{]}

Let $\omega\models q,\alpha\models p$.

Note that since $p$ is not definable over $M$, for every $m\in M$,
and even for every $m\in\C$ such that $\tp\left(m/M\right)$ is definable,
there is some $\varepsilon_{m}\in M$ such that $0<\varepsilon_{m}<\left|\alpha-m\right|$. 

Now, suppose that $\varphi\left(x,y\right)$ is any formula over $M$.
Then, as $q$ is definable, there is a formula $\psi\left(x\right):=\left(d_{q}y\right)\varphi\left(x,y\right)$
over $M$ that defines $\varphi\left(M,\omega\right)$. We claim that
$p\left(x\right)\cup q\left(y\right)\models\varphi\left(x,y\right)$
if and only if $p\left(x\right)\models\left(d_{q}y\right)\varphi\left(x,y\right)$.

We know that $\psi$ is equivalent to a finite union of intervals
and points from $M$. We also know that $\varphi\left(\C,\omega\right)$
is such a union, but the types of the end-points over $M$ are definable
over $M$ (since we have definable Skolem functions). So denote the
set of all these end-points by $A$. Let $0<\varepsilon\in M$ be
smaller than every $\varepsilon_{m}$ for each $m\in A$. Let $a,b\in M$
such that $a<\alpha<b$ and $b-a<\varepsilon$. Then:
\begin{itemize}
\item $\psi\left(\alpha\right)$ holds if and only if
\item $\psi\left(m\right)$ holds for all $m\in M$ such that $a\leq m\leq b$
if and only if
\item $\varphi\left(m,\omega\right)$ holds for all $m\in M$ such that
$a\leq m\leq b$ if and only if
\item $\varphi\left(\alpha,\omega\right)$ holds.
\end{itemize}
\end{proof}
We claim that RCF has large directionality. Moreover, we seem to answer
an open question raised in \cite{Delon} (at least in some sense,
see below), as she says there:
\begin{quotation}
Mais il laisse ouverte la possibilit\'e que la borne du nombre de
coh\'eritiers soit $\ded\left|M\right|$ dans le cas de la propri\'et\'e
de l'ordre et $\left(\ded\left|M\right|\right)^{\left(\omega\right)}$
dans le cas de l'ordre multiple.
\end{quotation}
So let us make clear what the question means:
\begin{defn}
(Taken from \cite{KeislerStabFunction,KeislerSixClasses}) A theory
$T$ is said to have the \emph{multiple order property} if there are
formulas $\varphi_{n}\left(x,y_{n}\right)$ for $n<\omega$ such that
the following set of formulas is consistent with $T$:
\[
\Gamma=\left\{ \varphi_{n}\left(x_{\eta},y_{n,k}\right)^{\eta\left(k\right)<n}\left|\,\eta\in\leftexp{\omega}{\omega}\right.\right\} .
\]
\end{defn}
\begin{rem}
If $T$ is strongly dependent (see \ref{def:StronglyDep}), for example,
if $T=RCF$, it does not have the multiple order property.\end{rem}
\begin{proof}
Suppose $T$ has the multiple order property as witnessed by formulas
$\varphi_{n}$. Consider the formulas $\psi_{n}\left(x,y,z\right)=\varphi_{n}\left(x,y\right)\leftrightarrow\varphi_{n}\left(x,z\right)$.
It is easy to see that $\left\{ \psi_{n}\left|\, n<\omega\right.\right\} $
exemplify that the theory is not strongly dependent. \end{proof}
\begin{fact}
\cite{KeislerStabFunction} If $T$ is countable and has the multiple
order property, then for every cardinal $\lambda$, $\sup\left\{ S\left(M\right)\left|\, M\models T,\left|\, M\right|=\lambda\right.\right\} \geq\left(\ded\lambda\right)^{\omega}$.
If $T$ does not have the multiple order property, then $\sup\left\{ S\left(M\right)\left|\, M\models T,\left|M\right|=\lambda\right.\right\} \leq\ded\lambda$. 
\end{fact}
So the question can be formulated as follows:
\begin{itemize}
\item Is there a countable theory without the multiple order property such
that for every $\lambda\geq\aleph_{0}$, $\sup\left\{ \left|\uf\left(p\right)\right|\left|\, p\in S_{<\omega}\left(M\right),\left|M\right|=\lambda\right.\right\} =\left(\ded\lambda\right)^{\omega}$
(recall that $S_{<\omega}\left(M\right)$ is the set of all finitary
types over $M$). 
\end{itemize}
It is a natural question, because of 2 reasons:
\begin{enumerate}
\item In general, the number of types (in $\alpha$ variables) over a model
of size $\lambda$ in a dependent theory is bounded by $\left(\ded\lambda\right)^{\left|T\right|+\left|\alpha\right|+\aleph_{0}}$
(by \cite[Theorem 4.3]{Sh10}), so by Fact \ref{fac:pOmega} this
is an upper bound for $\sup\left\{ \left|\uf\left(p\right)\right|\left|\, p\in S\left(M\right),\left|M\right|=\lambda\right.\right\} $.
\item It is very easy to construct an example with the multiple order property
that attains this maximum: for example, one can modify example \ref{exa:Large},
and add $\aleph_{0}$ independent orderings to $Q$.\end{enumerate}
\begin{defn}
For $M\models RCF$, let $S_{\dense}\left(M\right)$ be the set of
dense complete types over $M$.\end{defn}
\begin{thm}
\label{thm:more coheirs than dense types}Suppose $M\models RCF$.
Then there is a type $p\in S_{2}\left(M\right)$ such that $\left|\uf\left(p\right)\right|\geq\left|S_{\dense}\left(M\right)\right|^{\omega}$.\end{thm}
\begin{proof}
We may assume $S_{\dense}\left(M\right)\neq\emptyset$. Suppose $r_{*}$
is a dense type. Let $\alpha\models r_{*}$, and let $\omega\in\C$
be an element greater than any element in $M$. Then $q=\tp\left(\omega/M\right)$
is definable and we can apply Claim \ref{cla:WeakOrth}. Let $p(x_{\omega},x_{\alpha})=\tp(\omega,\alpha/M)$. 

For every dense type $r\left(x\right)$ over $M$, choose a realization
$a_{r}\in\C$. For every sequence $\bar{r}=\left\langle r_{i}\left|\, i<\omega\right.\right\rangle $
of \uline{positive} dense types over $M$ (i.e., $r_{i}\models x>0$)
, we define a coheir $p_{\bar{r}}$ of $p$ as follows: 

Fix $\bar{a}=\sequence{a_{r_{i}}}{i<\omega}$. For every sequence
$\bar{b}=\left\langle b_{i}\left|\, i<\omega\right.\right\rangle \in M^{\omega}$
such that $r_{i}\models x<b_{i}$ for all $i<\omega$, and for each
$n<\omega$ let $f_{n}\left(\bar{a},x\right)=\alpha+\sum_{i=0}^{n}\left(a_{r_{i}}/x^{i+1}\right)$
and $g_{n}\left(\bar{a},\bar{b},x\right)=\alpha+\sum_{i=0}^{n-1}\left(a_{r_{i}}/x^{i+1}\right)+b_{n}/x^{n+1}$.

Now, let $p_{\bar{r}}\left(x_{\omega},x_{\alpha}\right)$ be: 
\begin{eqnarray*}
p_{\bar{r}}\left(x_{\omega},x_{\alpha}\right) & = & p\left(x_{\omega},x_{\alpha}\right)\cup\left\{ f_{n}\left(\bar{a},x_{\omega}\right)<x_{\alpha}<g_{n}\left(\bar{a},\bar{b},x_{\omega}\right)\left|\,\bar{b}\mbox{ as above, }n<\omega\right.\right\} .
\end{eqnarray*}

\begin{claim*}
$p_{\bar{r}}\left(x_{\omega},x_{\alpha}\right)$ (which is over $M\cup\left\{ \alpha\right\} \cup\set{a_{r_{i}}}{i<\omega}$)
is finitely satisfiable in $M$. \end{claim*}
\begin{proof}
Suppose we are given a finite subset $p_{0}\subseteq p_{\bar{r}}\left(x_{\omega},x_{\alpha}\right)$,
and a finite set of inequalities $S=\set{f_{k}\left(\bar{a},x_{\omega}\right)<x_{\alpha}<g_{k}\left(\bar{a},\bar{b},x_{\omega}\right)}{k\leq n,\,\bar{b}\in B}$
where $B$ is some finite set of tuples $\sequence{b_{i}}{i\leq n}$
such that $r_{i}\models x<b_{i}$ for $i\leq n$. 

Let $\bar{b}$ be a tuple $\sequence{b_{i}}{i\leq n}\in M^{n+1}$
such that for $i\leq n$, $r_{i}\models x<b_{i}$ and $b_{i}<b_{i}'$
for any tuple $\sequence{b_{i}'}{i\leq n}\in B$. We may assume that
$p_{0}=r_{*,0}\left(x_{\alpha}\right)\cup q_{0}\left(x_{\omega}\right)$
where $r_{*,0}\subseteq r_{*}$ and $q_{0}\subseteq q$. We may assume
in addition that both $r_{*,0}$ and $q_{0}$ are intervals over $M$
(i.e., types in the language $\left\{ <\right\} $). Finally, we may
assume that $B=\left\{ \bar{b}\right\} $.

We will show:
\begin{itemize}
\item [$\smiley$]For all $o\in M$ large enough, there is some $0<\varepsilon_{o}\in M$
such that for all $k,l\leq n$, $\varepsilon_{o}<g_{k}\left(\bar{a},\bar{b},o\right)-f_{l}\left(\bar{a},o\right)$.
\end{itemize}
Once $\smiley$ is established, let $o$ be large enough so that it
has such an $\varepsilon_{o}$, $o$ satisfies $q_{0}\left(x_{\omega}\right)$
and for every $k\leq n$, $f_{k}\left(\bar{a},o\right),g_{k}\left(\bar{a},\bar{b},o\right)\models r_{*,0}\left(x_{\alpha}\right)$
(so also every element between $f_{k}$ and $g_{k}$). Suppose $l\leq n$
is such that $f_{l}\left(\bar{a},o\right)$ is maximal and $k\leq n$
is such that $g_{k}\left(\bar{a},\bar{b},o\right)$ is minimal. For
$i\leq n$, let $c_{i}\in M$ be such that $a_{r_{i}}<c_{i}$ and
$c_{i}-a_{r_{i}}<\varepsilon_{o}\cdot\left(o^{i+1}\right)/\left(l+2\right)$
(these exist since the $r_{i}$'s are dense), and let $\alpha<\alpha_{0}\in M$
be such that $\alpha_{0}-\alpha<\varepsilon_{o}/\left(l+2\right)$.
Let $d=\alpha_{0}+\sum_{i=0}^{l}\left(c_{i}/o^{i+1}\right)\in M$.
Then $f_{l}\left(\bar{a},o\right)<d$ and $d-f_{l}\left(\bar{a},o\right)=\left(\alpha_{0}-\alpha\right)+\sum_{i=0}^{l}\left(c_{i}-a_{r_{i}}\right)/o^{i+1}<\varepsilon_{o}$.
So $d<g_{k}\left(\bar{a},\bar{b},o\right)$, and so $\left(o,d\right)\models p_{0}$. 

So we only need to show $\smiley$. It is enough to show that for
each $k,l\leq n$, for all large enough $o$, there is some $0<\varepsilon_{o,k,l}\in M$
such that $\varepsilon_{o,k,l}<g_{k}\left(\bar{a},\bar{b},o\right)-f_{l}\left(\bar{a},o\right)$.
Suppose $k>l$. In that case, 
\[
g_{k}\left(\bar{a},\bar{b},o\right)-f_{l}\left(\bar{a},o\right)\geq b_{k}/o^{k+1}>0
\]
(since the types $r_{i}$ are positive). Suppose $k\leq l$. So,
\[
g_{k}\left(\bar{a},\bar{b},o\right)-f_{l}\left(\bar{a},o\right)=\left(b_{k}-a_{r_{k}}\right)/o^{k+1}-\left(\sum_{i=k+1}^{l}a_{r_{i}}/o^{i+1}\right).
\]
Since $r_{k}$ is dense, there is some $0<\varepsilon\in M$ such
that $\varepsilon<b_{k}-a_{r_{k}}$. Also, there are some $a_{i}'\in M$
such that $a_{r_{i}}<a_{i}'$. The difference above is greater than:
\[
\varepsilon/o^{k+1}-\sum_{i=k+1}^{l}a_{i}'/o^{i+1}\in M,
\]
and for $o$ large enough this number is positive, so let it be $\varepsilon_{o,k,l}$. 
\end{proof}
Note that for $\bar{r}\neq\bar{r}'$, $p_{\bar{r}}\cup p_{\bar{r}'}$
is inconsistent. Also, since $r_{*},r_{*}+1,r_{*}+2,\ldots$ are all
dense types, $\left|S_{\dense}\left(M\right)\right|\geq\aleph_{0}$,
so the number of positive dense types over $M$ is equal to the number
of all dense types over $M$. Together, we are done. 
\end{proof}
We conclude:
\begin{cor}
\label{thm:RCF-has-large}RCF has large directionality. In addition,
RCF does not have the multiple order property but for every $\lambda\geq\aleph_{0}$,
with $\cof\left(\ded\lambda\right)>\omega$, 
\[
\sup\left\{ \left|\uf\left(p\right)\right|\left|\, M\models RCF,\, p\in S_{2}\left(M\right),\,\left|M\right|=\lambda\right.\right\} \geq\left(\ded\lambda\right)^{\omega}.
\]
\end{cor}
\begin{proof}
We will use results from Section \ref{sec:Appendix:-dense-types}. 

By Theorem \ref{thm:ded via dense types}, we know that: 
\[
\sup\left\{ \left|\uf\left(p\right)\right|\left|\, p\in S_{2}\left(M\right),\left|M\right|=\lambda\right.\right\} \geq\sup\set{\left|S_{\dense}\left(M\right)\right|^{\omega}}{\left|M\right|=\lambda}.
\]

On the other hand, Corollary \ref{cor:what can we get with S_d} says
that: 
\[
\sup\set{\left|S_{\dense}\left(M\right)\right|^{\omega}}{\left|M\right|=\lambda}=\sup\set{\left(\lambda^{\left\langle \mu\right\rangle _{\tr}}\right)^{\omega}}{\mu\leq\lambda,\cof\left(\mu\right)=\mu},
\]
so this already implies that RCF is large (by Fact \ref{fac:calculating ded},
the right hand side is $\geq\ded\lambda$). Corollary \ref{cor:cof uncountable is ok}
says that for any cardinal $\lambda$, if $\cof\left(\ded\lambda\right)>\aleph_{0}$,
then 
\[
\sup\set{\left(\lambda^{\left\langle \mu\right\rangle _{\tr}}\right)^{\omega}}{\mu\leq\lambda,\cof\left(\mu\right)=\mu}=\left(\ded\lambda\right)^{\omega}.
\]
Together we are done. \end{proof}
\begin{rem}
For an easy proof that RCF is large, using the same notation from
the proof of Theorem \ref{thm:more coheirs than dense types}, for
every bounded cut $I\subseteq M$, define: 
\[
p_{I}\left(x_{\omega},x_{\alpha}\right)=r_{*}\left(x_{\alpha}\right)\cup q\left(x_{\omega}\right)\cup\set{\alpha+a/x_{\omega}<x_{\alpha}<\alpha+b/x_{\omega}}{a\in I,b\notin I}.
\]
 
\end{rem}
Marcus Tressl has pointed out the type $\tp\left(\alpha,\omega\right)$
to us as a type with infinitely many coheirs (this follows from \cite[18.13]{Tressle}).
We thank him for that. This proof that the theory is large is ours.

\subsubsection{Valued fields}

We can combine the techniques of Theorem \ref{thm:RCF-has-large}
and Example \ref{exa:Large} in order to prove a similar result for
valued fields. 
\begin{defn}
\label{def: language of val field}The language $L$ of valued fields
is the following. It is a 3-sorted language, one sort for the base
field $K$ equipped with the ring language $\left\{ 0,1,+,\cdot\right\} $,
another for the valuation group $\Gamma$ equipped with the ordered
abelian groups language $L_{\Gamma}=\left\{ 0,+,<\right\} $, and
another for the residue field $k$ equipped with the ring language
$L_{k}$. We also have the valuation map $v:K^{\times}\to\Gamma$
and an angular component map $\ac:K\to k$. Recall that an angular
component is a function that satisfies $\ac\left(0\right)=0$ and
$\ac\upharpoonright K^{\times}:K^{\times}\to k^{\times}$ is a homomorphism
such that if $v\left(x\right)=0$ then $\ac\left(x\right)$ is the
residue of $x$.
\end{defn}
For more on valued fields with angular component, see e.g., \cite{belair,Pas1989},
which also gives us the following fact:
\begin{fact}
\label{fac: Elimination of Field Quantifiers}\cite[Theorem 4.1]{Pas1989}
The theory of any Henselian valued field of characteristic $\left(0,0\right)$
in the language $L$ has elimination of field quantifiers: every formula
$\varphi\left(x_{K},x_{k},x_{\Gamma}\right)$ (where $x_{K}$, $x_{k}$
and $x_{\Gamma}$ are tuples of variables in the base field, the residue
field and the valuation group respectively) is equivalent to a Boolean
combination of formulas of the form $\varphi\left(\ac\left(f_{0}\left(x_{K}\right)\right),\ldots,\ac\left(f_{n-1}\left(x_{K}\right)\right),x_{k}\right)$
and $\chi\left(v\left(g_{0}\left(x_{K}\right)\right),\ldots,v\left(g_{m-1}\left(x_{K}\right)\right),x_{\Gamma}\right)$
where $\varphi$ is a formula in $L_{k}$, $\chi$ is a formula in
$L_{\Gamma}$ and $g_{i}$ and $f_{j}$ are polynomials over the integers. \end{fact}
\begin{thm}
Let $T=Th\left(K,\Gamma,k\right)$ be any theory of valued fields
in $L$ which eliminates field quantifiers. Then $T$ has large directionality.\end{thm}
\begin{proof}
Let $M_{0}=\left(K_{0},\Gamma_{0},k_{0}\right)\models T$ be a countable
model such that $\Gamma_{0}$ contains a copy of the rationals $\set{\gamma_{q}}{q\in\Qq}$
with the usual order and group structure, so $\gamma_{0}=0_{\Gamma}$
(by compactness, one only needs to embed a copy of a finitely generated
subgroup of $\left(\Qq,+,<\right)$ in a model of $T$, but any such
subgroup is contained in a subgroup generated by one element, which
is isomorphic to $\left(\mathbb{Z},+,<\right)$).

Let $S$ be the tree $2^{\leq\omega}$ and let $S_{0}\subseteq S$
be $2^{<\omega}$, so $S_{0}$ is countable. 

Let $\Sigma\sequence{x_{s}}{s\in S_{0}}$ be the following set of
formulas with variables in the field sort (over $M_{0}$): 
\[
\left\{ v\left(x_{s}-x_{t}\right)=\gamma_{\lev\left(s\wedge t\right)}\left|\, s,t\in S_{0},\, s\wedge t<s,t\right.\right\} \cup\set{v\left(x_{s}-x_{t}\right)\geq\gamma_{\lev\left(s\right)}}{s,t\in S_{0},\, s\leq t}.
\]
Then $\Sigma$ is consistent with $M_{0}$: to realize $\Sigma\upharpoonright2^{<n}$,
choose $a_{i}\in K_{0}$ with $v\left(a_{i}\right)=\gamma_{i}$ and
let $x_{s}=\sum_{i<\lev\left(s\right)}s\left(i\right)a_{i}$ for $s\in2^{<n}$.
Let $M=\left(K_{1},\Gamma_{1},k_{1}\right)$ be a countable model
containing $M_{0}$ and some $\set{a_{s}}{s\in S_{0}}$ realizing
$\Sigma$. 

For each $\eta\in S$ with domain $\omega$ (this is a \emph{branch}
of $S_{0}$), let $p_{\eta}\left(x\right)$ be the following type
in the valued field sort: 
\[
\left\{ v\left(x-a_{s}\right)\geq\gamma_{\lev\left(s\right)}\left|\, s<\eta\right.\right\} .
\]
 It is consistent since any finite subset if realized by $a_{t}$
for any $t<\eta$ large enough. 

If $\eta_{1}\neq\eta_{2}$ then $p_{\eta_{1}}\cup p_{\eta_{2}}$ is
inconsistent: 

Suppose $s=\eta_{1}\wedge\eta_{2}$ and $s<t<\eta_{1}$, $s<t'<\eta_{2}$.
If $p_{\eta_{2}}$ is consistent with $v\left(x-a_{t}\right)\geq\gamma_{\lev\left(t\right)}$,
then there is some $a$ such that $v\left(a-a_{t}\right)\geq\gamma_{\lev\left(t\right)}$
and $v\left(a-a_{t'}\right)\geq\gamma_{\lev\left(t'\right)}$. So
$v\left(a_{t}-a_{t'}\right)\geq\min\left\{ \gamma_{\lev\left(t'\right)},\gamma_{\lev\left(t\right)}\right\} $,
but $t\wedge t'=s<t,t'$ and so $v\left(a_{t}-a_{t'}\right)=\gamma_{\lev\left(s\right)}$.
This is a contradiction since $\gamma_{\lev\left(t\right)},\gamma_{\lev\left(t'\right)}>\gamma_{\lev\left(s\right)}$.

Let $\Omega$ be the algebraic closure (as a valued field) of the
monster model $\C$ of $M$. let $\bar{M}=\left(\bar{K_{1}},\bar{\Gamma_{1}},\bar{k_{1}}\right)$
be the algebraic closure of $M$ as a valued field in $\Omega$. Since
$\bar{M}$ is countable, there is some branch $\eta\in S$ such that
$p_{\eta}$ is not realized in $\bar{M}$. Then for every polynomial
$f\left(x\right)$ over $K_{1}$ and every $s<\eta$ large enough,
$p_{\eta}\models\ac\left(f\left(x\right)\right)=\ac\left(f\left(a_{s}\right)\right)\land v\left(f\left(x\right)\right)=v\left(f\left(a_{s}\right)\right)$
(decompose $f$ into linear factors $\prod\left(x-b_{i}\right)$.
For every large enough $s$, $v\left(b_{i}-a_{s}\right)<\gamma_{\lev\left(s\right)}$
for all $i$, and so if $d\models p_{\eta}$ in $\C$, then $\ac\left(f\left(d\right)\right)=\ac\left(f\left(a_{s}\right)\right)$
(because $res\left(\frac{f\left(d\right)}{f\left(a_{s}\right)}\right)=1$
--- we do not assume that $\ac$ extends to $\bar{M}$) and $v\left(f\left(d\right)\right)=v\left(f\left(a_{s}\right)\right)$).
Since field quantifiers are eliminated in $T$, this implies that
$p_{\eta}$ is a complete type. Moreover, we have the following claim:
\begin{claim*}
For any type $r\left(y\right)\in S\left(M\right)$ such that $y$
is a tuple of variables in the valuation group sort, $p_{\eta}$ and
$q$ are weakly orthogonal, meaning that $p_{\eta}\left(x\right)\cup r\left(y\right)$
implies a complete type over $M$.\end{claim*}
\begin{proof}
By elimination of field quantifiers, we need only to determine whether
\[
\chi\left(v\left(g_{0}\left(x\right)\right),\ldots,v\left(g_{m-1}\left(x\right)\right),y\right)
\]
 is in $p_{\eta}\left(x\right)\cup q\left(y\right)$ for any formula
$\chi$ in $L_{\Gamma}$ over $M$ and polynomials $g_{i}$ over $M$.
By the remark above, $p_{\eta}\models v\left(g_{i}\left(x\right)\right)=v\left(g_{i}\left(a_{s}\right)\right)$
for any $s<\eta$ large enough and all $i<m$. So $p_{\eta}\cup r\models\chi$
iff $r\models\chi\left(v\left(g_{0}\left(a_{s}\right)\right),\ldots,v\left(g_{m-1}\left(a_{s}\right)\right),y\right)$. 
\end{proof}
Let $r\left(y\right)\in S\left(M\right)$ be a type in the valuation
group sort which is finitely satisfiable in $\set{\gamma_{q}}{q\in\Qq}$
and contains $\set{y>\gamma_{q}}{q\in\Qq}$. By the claim, $r$ and
$p_{\eta}$ are weakly orthogonal. Fix some $d\models p_{\eta}$ in
$\C$. For each bounded cut $I\subseteq\Qq$, let $p_{I}\left(x,y\right)$
be the following type:
\[
p_{\eta}\left(x\right)\cup r\left(y\right)\cup\set{y+\gamma_{q}<v\left(x-d\right)<y+\gamma_{q'}}{q\in I,\, q'\notin I}.
\]

Then $p_{I}\left(x,y\right)$ is finitely satisfiable in $M$: 

Suppose we are given finite subsets $p_{0}\subseteq p_{\eta}$ and
$r_{0}\subseteq r$, $I_{0}\subseteq I$ and $I_{0}'\subseteq\Qq\backslash I$.
Let $q=\max I_{0}$ and $q'=\min I_{0}'$. Note that there is some
$s<\eta$ such that for any $a\in K_{1}$, if $v\left(a-a_{s}\right)\geq\gamma_{\lev\left(s\right)}$
then $a\models p_{0}$. Let $q_{0}\in\Qq$ be larger than $\lev\left(s\right)$,
larger than $\lev\left(s\right)-q$ and such that $\gamma_{q_{0}}\models r_{0}$.
Let $q''\in\Qq$ be in the interval $\left(q_{0}+q,q_{0}+q'\right)$
and let $b\in K_{1}$ be such that $v\left(b\right)=\gamma_{q''}$.
Let $s<t<\eta$ be such that $\lev\left(t\right)>q''$ and let $a=a_{t}+b\in K_{1}$.
Then $p_{0}\left(a\right)\cup r_{0}\left(\gamma_{q_{0}}\right)$ holds,
and in addition, 
\[
v\left(a-d\right)=v\left(a_{t}-d+b\right)=v\left(b\right)=\gamma_{q''}.
\]
Moreover,
\[
\gamma_{q_{0}}+\gamma_{q}=\gamma_{q_{0}+q}<\gamma_{q''}<\gamma_{q_{0}+q'}=\gamma_{q_{0}}+\gamma_{q'}.
\]

Obviously, for different cuts $I$ and $J$, the types $p_{I}$ and
$p_{J}$ contradict each other. 

Together this shows that, letting $p\left(x,y\right)$ be the complete
type determined by $p_{\eta}\left(x\right)\cup r\left(y\right)$,
$\left|\uf_{\Delta}\left(p\right)\right|\geq2^{\aleph_{0}}$ where
$\Delta=\left\{ \varphi\left(x,y;z_{0},z_{1},z_{2}\right)\right\} $
and: 
\begin{eqnarray*}
\varphi\left(x,y;z_{0},z_{1},z_{2}\right) & = & y+z_{1}<v\left(x-z_{0}\right)<y+z_{2}.
\end{eqnarray*}
So $T$ is large as promised. 
\end{proof}
There are other languages of valued fields in addition to the one
in Definition \ref{def: language of val field} that would make the
proof above work. The only requirements are that we can construct
the tree $S$ inside the field, that $\Gamma$ is a sort and that
field quantifiers are eliminated. This can be done in the theory of
the $p$-adics, when we add to the language $\ac_{n}$ for $n<\omega$
as in \cite{MR1023804}, and also in ACVF --- algebraically closed
valued fields (where there is quantifier elimination in any reasonable
language, and in fact there is no need for $\ac$). 
\begin{cor}
The theory of any Henselian valued field of characteristic $\left(0,0\right)$
(in the language described in Definition \ref{def: language of val field}),
ACVF (in any characteristic and any reasonable language with quantifier
elimination and a sort for the valuation group), and the theory of
$\Qq_{p}$ (in the language of Pas with $\ac_{n}$) are large. 
\end{cor}

\section{\label{sec:Splintering}Splintering}

This part of the paper is motivated by the work of Rami Grossberg,
Andr\'es Villaveces and Monica VanDieren. In their paper \cite{GrViVa}
they study Shelah's Generic pair conjecture (which is now a theorem
--- \cite{Sh:900,Sh950,Sh906}), and in their analysis they came up
with the notion of splintering, a variant of splitting.
\begin{defn}
Let $p\in S\left(\C\right)$. Say that $p$ \emph{splinters} over
$M$ if there is some $\sigma\in\Aut\left(\C\right)$ such that 
\begin{enumerate}
\item $\sigma\left(p\right)\neq p$.
\item $\sigma\left(p|_{M}\right)=p|_{M}$.
\item $\sigma\left(M\right)=M$ setwise.
\end{enumerate}
\end{defn}
\begin{rem}
{[}due to Martin Hils{]} Splitting implies splintering, and if $T$
is stable, then they are equal. \end{rem}
\begin{proof}
Suppose $p\in S\left(\C\right)$ does not split over $M$, then, by
stability, it is definable over $M$, and $p$ is the unique non-forking
extension of $p|_{M}$. Then for any $\sigma\in\Aut\left(\C\right)$,
$\sigma\left(p\right)$ is the unique non-forking extension of $\sigma\left(p|_{M}\right)$.
So if $\sigma\left(p|_{M}\right)=p|_{M}$, this means that $\sigma\left(p\right)=p$
so $p$ does not splinter over $M$.\end{proof}
\begin{claim}
\label{cla:unstable implies sp neq spln}Outside of the stable context,
splitting $\neq$ splintering.\end{claim}
\begin{proof}
Let $T$ be the theory of random graphs in the language $\left\{ I\right\} $.
Let $M\models T$ be countable, and let $a\neq b\in M$ with an automorphism
$\sigma\in\Aut\left(M\right)$ taking $a$ to $b$. Let $p\left(x\right)\in S\left(\C\right)$
say that $x\mathrela Ic$ for every $c\in M$ and if $c\notin M$
then $x\mathrela Ic$ if and only if $c$ is connected $a$ and not
connected to $b$. Obviously, $p$ does not split over $M$. However,
let $\sigma'\in\Aut\left(\C\right)$ be an extension of $\sigma$.
Let $c\in\C$ be such that $c$ is connected to $a$ but not to $b$.
Then $x\mathrela Ic\in p$ but $x\mathrela Ic\notin\sigma'\left(p\right)$. 
\end{proof}
However,
\begin{claim}
If $T=Th\left(\mathbb{Q},<\right)$, then splitting equals splintering.\end{claim}
\begin{proof}
Observe that by quantifier elimination every complete type $r\left(x_{i}\left|\, i\in I\right.\right)$
over a set $A$ is determined by $\tp\left(\sequence{x_{i}}{i\in I}/\emptyset\right)\cup\bigcup\left\{ \tp\left(x_{i}/A\right)\left|\, i\in I\right.\right\} $.
Assume $q\left(x_{i}\left|\, i\in I\right.\right)$ is a global type
that splinters but does not split over a model $M$. Then it follows
that for some $i\in I$, $q\upharpoonright x_{i}$ splinters, so we
may assume $\left|I\right|=1$. Suppose $\sigma\in\Aut\left(\C\right)$
is such that $\sigma\left(M\right)=M$, $\sigma\left(q|_{M}\right)=q|_{M}$
and $\sigma\left(q\right)\neq q$. Note that $\sigma\left(q^{\left(\omega\right)}\right)=\sigma\left(q\right)^{\left(\omega\right)}$,
so by Fact \ref{fac:pOmega}, $\sigma\left(q^{\left(\omega\right)}\right)|_{M}\neq q^{\left(\omega\right)}|_{M}$.
We get a contradiction by quantifier elimination again. 
\end{proof}
We shall now generalize Claim \ref{cla:unstable implies sp neq spln}
to every theory with the independence property. In fact, to any theory
with large or medium directionality.
\begin{thm}
Suppose $T$ has medium or large directionality then splitting $\neq$
splintering.\end{thm}
\begin{proof}
We know that there is some $p$ and $\Delta$ such that $\uf_{\Delta}\left(p\right)$
is infinite. Let us use Construction \ref{const:not small directionality}:

We may find a saturated model $N=\left(N_{0}',N_{0},M_{0},Q_{0},\bar{f}_{0}\right)$
of $Th\left(M^{*}\right)$ of size $\lambda$ where $\lambda$ is
big enough. Then there is $c\neq d\in Q_{0}$ such that $\tp\left(c/\emptyset\right)=\tp\left(d/\emptyset\right)$
in the extended language (with symbols for $N_{0},M_{0},Q_{0}$ and
$\bar{f}_{0}$). So there is an automorphism $\sigma$ of this structure
(in particular of $N_{0}'$) such that $\sigma\left(c\right)=d$.
By definition, $\sigma\left(N_{0}\right)=N_{0}$ and $\sigma\left(M_{0}\right)=M_{0}$.
So $\tp\left(c/N_{0}\right)$ is finitely satisfiable in $M_{0}$
and hence does not split over $M_{0}$. But it splinters since $\sigma\left(\tp\left(c/M_{0}\right)\right)=\tp\left(d/M_{0}\right)=\tp\left(c/M_{0}\right)$
but $\sigma\left(\tp\left(c/N_{0}\right)\right)=\tp\left(d/N_{0}\right)\neq\tp\left(c/N_{0}\right)$
as witnessed by $\varphi$.

If there are no saturated models, we can take a big enough special
model (see \cite[Theorem 10.4.4]{Hod}). 

Note that we may also find an example of a type $p\in S\left(M\right)$
with a splintering, non-splitting, global extension, with $\left|M\right|=\left|T\right|$:
consider the structure $\left(N_{0}',N_{0},M_{0},\sigma,c,d\right)$,
and find an elementary substructure of size $\left|T\right|$. \end{proof}
\begin{defn}
Let $T$ be a complete theory. We say that $\left(M,p,\varphi\left(x;y\right),A_{1},A_{2}\right)$
is an\emph{ sp-example} for $T$ when:
\begin{itemize}
\item $M\models T$; $A_{1},A_{2}\subseteq M$ are nonempty and disjoint;
$p=p\left(x\right)$ is a complete type over $M$, finitely satisfiable
in $A_{1}$; $\Th\left(M_{p},A_{1}\right)=\Th\left(M_{p},A_{2}\right)$
(see Definition \ref{def:ExtDefType}); For each pair of finite sets
$s_{1}\subseteq A_{1}$ and $s_{2}\subseteq A_{2}$, $M\models\exists y\left(\bigwedge_{a\in s_{1}}\varphi\left(a,y\right)\land\bigwedge_{b\in s_{2}}\neg\varphi\left(b,y\right)\right)$.
\end{itemize}
\end{defn}
\begin{prop}
\label{pro:sp-example}$T$ has an sp-example if and only if there
is a finitely satisfiable type over a model which splinters over it
(in particular, splitting is different than splitting).\end{prop}
\begin{proof}
Suppose $\left(M,p,\varphi\left(x,y\right),A_{1},A_{2}\right)$ is
an sp-example for $T$. Let $M'$ be the structure $\left(M_{p},A_{1},A_{2}\right)$
(in the language $L_{p}\cup\left\{ P_{1},P_{2}\right\} $ where $P_{1},P_{2}$
are predicates). Assume $\left|T\right|<\mu=\mu^{<\mu}$, and let
$N'=\left(N_{q},B_{1},B_{2}\right)$ be a saturated extension of $M'$
of size $\mu$ where $N=N'\upharpoonright L$ and $q=q^{N}$ is as
in Remark \ref{rem:TpModel}. Since $\left(N_{q},B_{1}\right)\equiv\left(N_{q},B_{2}\right)$,
there is an automorphism $\sigma$ of $N_{q}$, such that $\sigma$
takes $B_{1}$ to $B_{2}$ and so $\sigma\left(q\right)=q$. Let $q'$
be a global extension of $q$, finitely satisfiable in $B_{1}$ and
$\sigma'$ a global extension of $\sigma$. 

So $q'$ does not split over $N$, but it splinters:

Consider the type $\left\{ \varphi\left(a,y\right)\left|\, a\in B_{1}\right.\right\} \cup\left\{ \neg\varphi\left(b,y\right)\left|\, b\in B_{2}\right.\right\} $.
It is finitely satisfiable in $N$ by choice of $\varphi$. Let $c\in\C$
satisfy this type. Then $\varphi\left(x,c\right)\in q'$ but $\varphi\left(x,c\right)\notin\sigma\left(q'\right)$
(because $\sigma\left(q'\right)$ is finitely satisfiable in $B_{2}$). 

If we do not assume the existence of such a $\mu$, we can use special
models.

Now suppose that splitting is different than splintering, as witnessed
by some global type $p$ that splinters over a model $M$ but is finitely
satisfiable in it. Then there is some automorphism $\sigma$ of $\C$
that witnesses it. There is a formula $\varphi\left(x,y\right)$ and
$a\in\C$ such that $\varphi\left(x,a\right)\in p,\neg\varphi\left(x,a\right)\in\sigma\left(p\right)$.
Let $B_{1}=\left\{ m\in M\left|\,\varphi\left(m,a\right)\land\neg\varphi\left(m,\sigma^{-1}\left(a\right)\right)\right.\right\} $,
$B_{2}=\sigma\left(B_{1}\right)$. It is easy to check that $\left(M,p,\varphi,B_{1},B_{2}\right)$
is an sp-example
\end{proof}
The following theorem answers the natural question:
\begin{thm}
There is a theory with small directionality in which splitting $\neq$
splintering.\end{thm}
\begin{proof}
Let $L=\left\{ R\right\} $ where $R$ is a ternary relation symbol.
Let $M_{0}=\left\langle \mathbb{Q},<\right\rangle $ and define $R\left(x,y,z\right)$
by $x<y<z$ or $z<y<x$, i.e., $y$ is between $x$ and $z$. Let
$T=Th\left(M_{0}\upharpoonright L\right)$. 
\begin{claim*}
$T$ has small directionality.\end{claim*}
\begin{proof}
Suppose $M\models T$. Let $\left(a,b\right)$ denote $\left\{ c\left|\, R\left(a,c,b\right)\right.\right\} $.

Then, for any choice of a pair of distinct elements $a,b$ there is
a unique enrichment of $M$ to a model $M'$ of $Th\left(\mathbb{Q},<\right)$
such that $R$ is defined as above and $a<b$: 

For $w\neq z$, $w<z$ if and only if

$\left(a,b\right)\cap\left(a,z\right)\neq\emptyset$ and ($\left(a,w\right)\subseteq\left(a,z\right)$
or $\left(a,w\right)\cap\left(a,z\right)=\emptyset$) or

$\left(a,b\right)\cap\left(a,z\right)=\emptyset$ and $\left(z,b\right)\subseteq\left(w,b\right)$.

From this observation, it follows that there is a unique completion
of any type $p\in S\left(M\right)$ to a type $p'\in S\left(M'\right)$.
So if $\Delta$ is a finite set of $L$ formulas and $\uf_{\Delta}\left(p\right)$
is infinite, then $\uf_{\Delta}\left(p'\right)$ is also infinite
--- contradiction to Example \ref{exa:DLO}.\end{proof}
\begin{claim*}
$T$ has an sp-example\end{claim*}
\begin{proof}
Let $M=M_{0}\upharpoonright L$. Let $p\left(x\right)=\tp\left(\pi/M\right)$.
Let $A_{1}=\left\{ x\in\mathbb{Q}\left|\, x>\pi\right.\right\} $
and $A_{2}=\mathbb{Q}\backslash A_{1}$, and let $\varphi\left(x;y_{1},y_{2}\right)=R\left(y_{1},x,y_{2}\right)$.
We claim that $\left(M,p,\varphi\left(x,y_{1},y_{2}\right),A_{1},A_{2}\right)$
is an sp-example:

First, let $M'$ be the reduct of $\left(\mathbb{Q}\cup\left\{ \pi\right\} ,<\right)$
to $L$. There is some $\sigma\in\Aut\left(M'/\pi\right)$ such that
$\sigma\left(A_{1}\right)=A_{2}$. Hence $\left(M_{p},A_{1}\right)\cong\left(M_{p},A_{2}\right)$.
Also, since $\tp\left(\pi/M_{0}\right)$ (in $\left\{ <\right\} $)
is finitely satisfiable in both $A_{1}$ and $A_{2}$ (by quantifier
elimination), $p$ is finitely satisfiable in both $A_{1}$ and $A_{2}$.
Finally, for finite $s_{1}\subseteq A_{1}$ and $s_{2}\subseteq A_{2}$,
there exists $c_{1},c_{2}\in\mathbb{Q}$ such that $R\left(c_{1},a,c_{2}\right)$
for all $a\in s_{1}$ and $\neg R\left(c_{1},b,c_{2}\right)$ for
all $b\in s_{2}$. 
\end{proof}
\end{proof}

\section{\label{sec:Appendix:-dense-types}Appendix: dense types in RCF}
\begin{defn}
For $M\models RCF$, let $S_{\dense}\left(M\right)$ be the set of
dense complete types over $M$ (see Definition \ref{def:dense types}). 
\end{defn}
Here we will prove the following theorem:
\begin{thm}
\label{thm:ded via dense types}$\ded\lambda=\sup\set{\left|S_{\dense}\left(M\right)\right|}{M\models RCF,\left|M\right|=\lambda}$.
\end{thm}
For the proof we will need some definitions and facts: 
\begin{defn}
\label{def:trees}
\begin{enumerate}
\item By a \emph{tree} we mean a partial order $\left(T,<\right)$ such
that for every $a\in T$, $T_{<a}=\left\{ x\in T\left|\, x<a\right.\right\} $
is well ordered. For $a\in T$, the order type of $T_{<a}$ is $a$'s
\emph{level}. By a \emph{branch} in $T$ we mean a maximally linearly
ordered subset of $T$. Its length is its order type. 
\item For two cardinals $\lambda$ and $\mu$, let $\lambda^{\left\langle \mu\right\rangle _{\tr}}$
be: 
\[
\sup\left\{ \kappa\left|\,\mbox{\mbox{there is some tree} }T\mbox{ with }\lambda\mbox{ many nodes and }\kappa\mbox{ branches of length }\mu\right.\right\} .
\]

\end{enumerate}
\end{defn}
\begin{fact}
\label{fac:calculating ded}(See \cite[Theorem 2.1(a)]{Baumgartner})
The following cardinalities are the same: 
\begin{enumerate}
\item $\ded\lambda$.
\item $\sup\left\{ \kappa\left|\,\mbox{\mbox{there is a regular }\ensuremath{\mu}\ and a tree }T\mbox{ with }\kappa\mbox{ branches of length }\mu\mbox{ and }\left|T\right|\leq\lambda\right.\right\} $.
\item $\sup\left\{ \lambda^{\left\langle \mu\right\rangle _{\tr}}\left|\,\mu\leq\lambda\mbox{ is regular}\right.\right\} $.
\end{enumerate}
\end{fact}
It is somewhat easier to consider trees which are sub-trees of $\lambda^{<\mu}$
(with the usual ``first-segment'' order) for some $\lambda,\mu$.
Given any tree $T$, and any cardinal $\mu$, suppose we are interested
in computing the number of branches of length $\mu$. For this we
may assume that the level of each element in $T$ is $<\mu$. Suppose
$\left|T\right|=\lambda$, so we may assume that its universe is $\lambda$.
Let $T'$ be $\set{T_{<a}}{a\in T}$. This is easily seen to be a
tree with the inclusion ordering, and moreover it is isomorphic to
a complete sub-tree $T''$ of $\lambda^{<\mu}$ (in the sense that
if $\eta\in T''$ and $\nu$ is an initial segment of $\eta$, then
$\nu\in T''$): if $\lev\left(a\right)=\alpha$, map $T_{<a}$ to
$\eta:\alpha\to\lambda$ where $\eta\left(\beta\right)$ is the $\beta$'th
element in $T_{<a}$. If $B\subseteq T$ is a branch of length $\mu$,
let $B'=\set{T_{<a}}{a\in B}$. Then $B'$ is also a branch of length
$\mu$, and in addition if $B_{1}\neq B_{2}$ are branches of $T$,
then $B_{1}'\neq B_{2}'$ in $T'$. This shows that $T'$ (so also
$T''$) has at least as many branches as $T$, and so in calculating
$\ded\lambda$ we can add to our list of cardinalities from Fact \ref{fac:calculating ded}:
\begin{enumerate}
\item [(4)]$\sup\left\{ \kappa\left|\,\mbox{\mbox{there is a regular }\ensuremath{\mu}\ and a tree }T\subseteq\lambda^{<\mu}\mbox{ with }\kappa\mbox{ branches of length }\mu\mbox{ and }\left|T\right|\leq\lambda\right.\right\} $. 
\end{enumerate}
Theorem \ref{thm:ded via dense types} follows from:
\begin{prop}
\label{prop:every tree manifests in dense types}For every tree $T\subseteq\lambda^{<\mu}$
of size $\lambda$, there is a model $M\models RCF$ of size $\lambda$
such that $\left|S_{\dense}\left(M\right)\right|$ is at least the
number of branches in $T$ of length $\mu$.\end{prop}
\begin{proof}
We may assume that $\mu\leq\lambda$. For $i<\mu$, let $T_{i}=T\cap\leftexp i{\lambda}$,
$T_{<i}=T\cap\lambda^{<i}$. By induction on $i<\mu$ we construct
a sequence of models $\bar{M}=\sequence{M_{i}}{i<\mu}$ and $\sequence{a_{\eta},b_{\eta}}{\eta\in T_{<i}}$
such that:

$\bar{M}$ is an $\prec$-increasing continuous sequence of models
of RCF; For all $\eta\in T_{<i}$, $a_{\eta},b_{\eta}\in M_{\lg\left(\eta\right)+1}$;
$a_{\eta}<b_{\eta}$; $b_{\eta}-a_{\eta}<c$ for all $0<c\in M_{\lg\left(\eta\right)}$;
If $\alpha<\beta<\lambda$ and $\eta\concat\left\langle \alpha\right\rangle ,\eta\concat\left\langle \beta\right\rangle \in T_{<i}$
then $b_{\eta\concat\left\langle \alpha\right\rangle }<a_{\eta\concat\left\langle \beta\right\rangle }$;
For $\nu<\eta$, $a_{\nu}<a_{\eta}<b_{\eta}<b_{\nu}$. 

The construction: 

Let $M_{0}$ be any model of size $\lambda$. 

For $i$ limit, let $M_{i}=\bigcup_{j<i}M_{j}$ (there are no new
$\left(a_{\eta},b_{\eta}\right)$'s). 

For $i=j+1$ for $j$ a successor, let $M_{i}$ be a model of size
$\lambda$ containing $M_{j}$ and an increasing sequence $\sequence{c_{\alpha}}{\alpha<\lambda}$
such that $0<c_{\alpha}<d$ for all $0<d\in M_{j}$. For $\eta\in T_{j-1}$,
if $\eta\concat\left\langle \alpha\right\rangle \in T$, let $a_{\eta\concat\left\langle \alpha\right\rangle }=a_{\eta}+c_{2\alpha}$
and $b_{\eta\concat\left\langle \alpha\right\rangle }=a_{\eta}+c_{2\alpha+1}$
(note that $b_{\eta\concat\left\langle \alpha\right\rangle }<b_{\eta}$). 

For $i=j+1$ for $j$ limit (or $j=0$), let $M_{i}$ be model of
size $\lambda$ containing $M_{j}$ and $a_{\eta},b_{\eta}$ for $\eta\in T_{j}$
where $a_{\eta\upharpoonright j'}<a_{\eta}<b_{\eta}<b_{\eta\upharpoonright j'}$
for all $j'<j$ (so for $j=0$ this just means $a_{\left\langle \right\rangle }<b_{\left\langle \right\rangle }$)
and $b_{\eta}-a_{\eta}<d$ for all $d\in M_{j}$. 

Finally, we let $M=\bigcup_{i<\mu}M_{i}$. For each branch $\eta\in\leftexp{\mu}{\lambda}$
of $T$, let $p_{\eta}=\set{a_{\eta\upharpoonright i}<x<b_{\eta\upharpoonright i}}{i<\mu}$.
This is easily seen to be a dense type. Also, it is very easy to see
that $p_{\eta}\neq p_{\eta'}$ for $\eta\neq\eta'$. \end{proof}
\begin{rem}
Note that this proof only used the fact that the order is dense, and
so this holds in any densely ordered abelian group. 
\end{rem}
Next we will show that Proposition \ref{prop:every tree manifests in dense types}
is ``as good as it gets''.
\begin{prop}
\label{prop:S_d is bounded}If $M\models RCF$, $\left|M\right|=\lambda$,
and $\mu=\cof\left(M,<\right)$, then $\left|S_{\dense}\left(M\right)\right|\leq\lambda^{\left\langle \mu\right\rangle _{\tr}}$.\end{prop}
\begin{proof}
We shall construct a tree of size $\lambda$ with $\left|S_{\dense}\left(M\right)\right|$
branches of length $\mu$. 

Let $\sequence{d_{i}}{i<\mu}$ be an increasing cofinal sequence of
positive elements in $M$. Let $<^{*}$ be a well ordering on $M^{2}$.
We define a sequence of pairs $\sequence{\left(a_{i,p},b_{i,p}\right)}{i<\mu,p\in S_{\dense}\left(M\right)}$
by induction on $i<\mu$ such that:

$\left(a_{i,p},b_{i,p}\right)\in M^{2}$ is the $<^{*}$-first pair
such that $p\left(x\right)\models a_{i,p}<x<b_{i,p}$, $b_{i,p}-a_{i,p}<1/d_{i}$
and for $j<i$, $a_{j,p}<a_{i,p}$, $b_{i,p}<b_{j,p}$. 
\begin{claim*}
$\left(a_{i,p},b_{i,p}\right)$ exist for all $p\in S_{\dense}\left(M\right)$
and $i<\mu$.\end{claim*}
\begin{proof}
Fix some $p\in S_{\dense}\left(M\right)$. Suppose $i<\mu$ is the
first such that $\left(a_{i,p},b_{i,p}\right)$ do not exist. For
$j<i$, let $0<c_{j}\in M$ be such that $p\left(x\right)\models x+c_{j}<b_{j,p}$
and $p\left(x\right)\models a_{j,p}+c_{j}<x$ (exists since $p$ is
dense, since otherwise it would be definable). Since the cofinality
of $M$ is $\mu$, there must be some $e\in M$ such that $e>d_{i}$
and $e>1/c_{j}$ for all $j<i$. Since $p$ is dense there must be
some $a,b\in M$ such that $p\left(x\right)\models a<x<b$ and $b-a<1/e$.
By choice of $e$ for all $j<i$, $a_{j,p}<a$, $b<b_{j,p}$ and $b-a<1/d_{i}$. 
\end{proof}
For $i<\mu$, let: 
\[
T_{i}=\set{\eta:i\to M^{2}}{\exists p\in S_{\dense}\left(M\right)\forall j<i\left[\eta\left(j\right)=\left(a_{j,p},b_{j,p}\right)\right]}.
\]

\begin{claim*}
If $\eta\in T_{i}$ then $\eta\upharpoonright j\in T_{j}$ for all
$j<i$. 
\end{claim*}

\begin{claim*}
$\left|T_{i}\right|\leq\lambda$. \end{claim*}
\begin{proof}
By the first claim, if $\eta\in T_{i}$ then it can be extended to
some $\nu$ in $T_{i+1}$. So it is enough to show that $\left|T_{i+1}\right|\leq\lambda$.
For that it is enough to show that the map $\eta\mapsto\eta\left(i\right)$
from $T_{i+1}$ to $M^{2}$ is injective. But this follows from definition
of $\left(a_{i,p},b_{i,p}\right)$. 
\end{proof}
Let $T=\bigcup_{i<\mu}T_{i}$. Then $T$ a tree, and for each dense
type $p\in S_{\dense}\left(M\right)$, we can find a branch $\eta_{p}:\mu\to M^{2}$
defined by $\eta\left(i\right)=\left(a_{i,p},b_{i,p}\right)$. The
following claim finishes the proof:
\begin{claim*}
For $p_{1}\neq p_{2}$, $\eta_{p_{1}}\neq\eta_{p_{2}}$.\end{claim*}
\begin{proof}
Suppose $p_{1}\left(x\right)\models x<b$ and $p_{2}\left(x\right)\models b<x$,
and let $0<e\in M$ be such that $p_{1}\left(x\right)\models x+e<b$
and $p_{2}\left(x\right)\models b+e<x$ (exists since $p_{1}$ and
$p_{2}$ are not definable). For some $i<\mu$, $d_{i}>1/e$. Then
it follows that $\eta_{p_{1}}\left(i\right)\neq\eta_{p_{2}}\left(i\right)$
. 
\end{proof}
\end{proof}
\begin{cor}
\label{cor:what can we get with S_d}The following equality holds
for all cardinals $\lambda\geq\aleph_{0}$:
\[
\sup\set{\left|S_{\dense}\left(M\right)^{\omega}\right|}{M\models RCF,\,\left|M\right|=\lambda}=\sup\set{\left(\lambda^{\left\langle \mu\right\rangle _{\tr}}\right)^{\omega}}{\mu\leq\lambda,\cof\left(\mu\right)=\mu}.
\]
\end{cor}
\begin{proof}
The inequality $\leq$ follows immediately from Proposition \ref{prop:S_d is bounded}.
For $\geq$ we will show that for every regular $\mu\leq\lambda$,
$\left(\treeexp{\lambda}{\mu}\right)^{\omega}\leq\sup\set{\left|S_{\dense}\left(M\right)^{\omega}\right|}{\left|M\right|=\lambda}$.

Suppose $\lambda^{\left\langle \mu\right\rangle _{\tr}}$ is attained,
i.e., there is a tree of size $\lambda$ with $\lambda^{\left\langle \mu\right\rangle _{\tr}}$
branches of length $\mu$. Then by Proposition \ref{prop:every tree manifests in dense types},
for some model $M\models RCF$ of size $\lambda$, $\left|S_{\dense}\left(M\right)\right|\geq\lambda^{\left\langle \mu\right\rangle _{\tr}}$,
so in that case we are done.

Suppose $\treeexp{\lambda}{\mu}$ is not attained. In that case $\cof\left(\treeexp{\lambda}{\mu}\right)>\lambda$.
Indeed, if not, then $\treeexp{\lambda}{\mu}=\bigcup_{i<\lambda}\sigma_{i}$
for some cardinals $\sigma_{i}<\treeexp{\lambda}{\mu}$. For each
$i<\lambda$, there is a tree $T_{i}$ of size $\lambda$ with more
than $\sigma_{i}$ branches of length $\mu$. Let $T$ be the disjoint
union of $T_{i}$ for $i<\lambda$. Then $T$ is a tree of size $\lambda$,
with at least $\treeexp{\lambda}{\mu}$ branches of length $\mu$
--- contradiction. In particular, $\cof\left(\treeexp{\lambda}{\mu}\right)>\omega$,
so every function $f:\omega\to\treeexp{\lambda}{\mu}$ is bounded,
and hence: 
\[
\left(\treeexp{\lambda}{\mu}\right)^{\omega}=\sup\set{\kappa^{\omega}}{\kappa<\treeexp{\lambda}{\mu}}.
\]
So it is enough to show that for each $\kappa<\treeexp{\lambda}{\mu}$,
there is a model $M$ of size $\lambda$ with more than $\kappa$
dense types, which follows from Proposition \ref{prop:every tree manifests in dense types}. \end{proof}
\begin{example}
In \cite[Section 6]{Sh1007} it is shown that it is consistent with
ZFC that there is an uncountable cardinal $\lambda$ such that:
\begin{enumerate}
\item $\cof\left(\ded\lambda\right)=\cof\left(\lambda\right)=\aleph_{0}$,
so $\left(\ded\lambda\right)^{\omega}>\ded\lambda$. 
\item For all regular cardinals $\mu<\lambda$, $\lambda^{\mu}\leq\ded\lambda$. 
\end{enumerate}
So in this case, 
\[
\sup\set{\left(\lambda^{\left\langle \mu\right\rangle _{\tr}}\right)^{\omega}}{\mu\leq\lambda,\cof\left(\mu\right)=\mu}=\ded\lambda.
\]

However, \end{example}
\begin{cor}
\label{cor:cof uncountable is ok}For any cardinal $\lambda$, if
$\cof\left(\ded\lambda\right)>\aleph_{0}$, then \textup{
\[
\sup\set{\left(\lambda^{\left\langle \mu\right\rangle _{\tr}}\right)^{\omega}}{\mu\leq\lambda,\cof\left(\mu\right)=\mu}=\left(\ded\lambda\right)^{\omega}.
\]
}\end{cor}
\begin{proof}
$\leq$ is clear. For $\geq$, we use a similar argument as in the
proof of Corollary \ref{cor:what can we get with S_d}. Since $\left(\ded\lambda\right)^{\omega}=\bigcup\set{\kappa^{\omega}}{\kappa<\ded\lambda}$,
we only have to show that every $\kappa<\ded\lambda$, $\kappa<\treeexp{\lambda}{\mu}$
for some regular $\mu\leq\lambda$. But that already follows from
Fact \ref{fac:calculating ded}. 
\end{proof}
\bibliographystyle{alpha}
\bibliography{common2}

\end{document}